\documentclass[10pt]{article}
\usepackage[latin1]{inputenc}
\usepackage{amsmath}
\usepackage{ epsfig}
\usepackage{graphicx}
\usepackage{epsf}
\usepackage{amsthm}
\usepackage{color}

\input epsf
\def\limitepsn{\renewcommand{\arraystretch}{0.5}
\begin{array}[t]{c}\stackrel{a.s.}{\longrightarrow} \\
{\scriptstyle n \rightarrow
\infty}\end{array}\renewcommand{\arraystretch}{1}}

\def\limiteloin{\renewcommand{\arraystretch}{0.5}
\begin{array}[t]{c}\stackrel{{\cal D}}{\longrightarrow} \\
{\scriptstyle n\rightarrow
\infty}\end{array}\renewcommand{\arraystretch}{1}}

\newcommand{\be}{\begin{equation}}
\newcommand{\ee}{\end{equation}}
\newcommand{\bd}{\begin{displaymath}}
\newcommand{\ed}{\end{displaymath}}

\newcommand{\ba}{\begin{eqnarray}}
\newcommand{\ea}{\end{eqnarray}}
\newcommand{\ban}{\begin{eqnarray*}}
\newcommand{\ean}{\end{eqnarray*}}

\newcommand{\R} {I\!\!R}
\newcommand{\E} {E\,}
\newcommand{\N} {I\!\! N}
\newcommand{\Z} {{\bf Z}}

\renewcommand{\arraystretch}{.8}

\renewcommand{\Box}{\hfill\rule{0.25cm}{0.25cm}} 

\newtheorem{lem}{Lemma}[section]
\newtheorem{Theo}{Theorem}[section]
\newtheorem{cor}{Corollary}[section]
\newtheorem{rem}{Remark}[section]
\theoremstyle{definition}
\newtheorem{exam}{Example}[section]
\newenvironment{dem}{\ \\ {\bf Proof. }}
{\Box\par\medskip\noindent}

\def\1{{\bf 1}}

\begin{document}
\title{\bf Monitoring procedure for parameter change in causal time series}
\author{Jean-Marc BARDET and William~KENGNE \footnote{Supported by Université Paris 1 and AUF (Agence Universitaire de la Francophonie).}}

 \maketitle \vspace{-1.0cm}
 %
%
\begin{center}
 {\it SAMM, Université Paris 1 Panth\'eon-Sorbonne, 90 rue de Tolbiac 75634-Paris Cedex 13, France.\\
  E-mail: Jean-Marc.Bardet@univ-paris1.fr ;  William-Charky.Kengne@malix.univ-paris1.fr
  }
\end{center}


 \textbf{Abstract} :   We propose a new sequential procedure to detect
    change in the parameters of a process $ X= (X_t)_{t\in \Z}$  belonging to a large class of causal models (such as AR($\infty$), ARCH($\infty$), TARCH($\infty$), ARMA-GARCH processes).
   The procedure is based on a difference between  the historical parameter estimator and the updated parameter estimator, where both these estimators are based on a quasi-likelihood of the model. Unlike classical recursive fluctuation test, the updated estimator is computed without the historical observations.
    The asymptotic behavior of the test is studied and the consistency in power as well as an upper bound of the detection delay are obtained.
    Some simulation results  are reported with comparisons to some other existing procedures exhibiting the accuracy of our new procedure.
    The procedure is also applied to the daily closing values of the Nikkei 225, S$\&$P 500 and FTSE 100 stock index. We show in this real-data applications
    how the procedure can be used to solve off-line multiple breaks detection.
~ \\

 {\em Keywords:}  Sequential change detection; Change-point; Causal
processes; Quasi-maximum likelihood estimator;  Weak convergence.

\section{Introduction} \label{intro}
In statistical inference, many authors have pointed out the danger of omitting the existence of changes in data.
Many papers have been devoted to the problem of test for parameter changes in time series models when all data are available, see for instance Horváth \cite{Horvart1993},
  Inclan and Tiao \cite{Tiao1994}, Kokoszka and Leipus \cite{Leipus1999}, Kim {\it et al.} \cite{Kim2000}, Horváth and Shao \cite{Horvath2007},
  Aue {\it et al.} \cite{Aue2009} or Kengne \cite{Kengne2011}. These papers consider "retrospective" (off-line)
  changes {\it i.e.} changes in parameters when all data are available. Another point of view is the change detection when new data arrive; this is the sequential change-point problem.
For instance, consider the following sequential problems.
  \begin{exam}[Industrial quality control]\label{ex1}
   Consider an industrial system producing electronic objects. To know in real time the quality of production, some devices have been installed to
   rely the informations about the functioning. The performance is evaluated for each production and the system is automatically stopped if a disorder is detected. After an investigation of the root cause of the
   problem and a possible maintenance, the system is resetting.
  \end{exam}
  \begin{exam}[Monitoring in economic system]\label{ex2}
    Consider an economic system described by the regression model $Y_t=X'_t\beta_0 + \xi_t, ~ t=1,2,\cdots n$ where $Y_t$ is the variable of interest and $X_t$ the vector of
    explanatory variables. New data arrive steadily at time $t=n+1,\cdots$; Chu {\it et al.} (1996) asked : "is yesterday's model capable of
    explaining today's data?" This question leads to construct a procedure which will test sequentially if the model with parameter $\beta_0$ is able to explain
    new data. If at time $k>n$ a change from $\beta_0$ to $\beta_1$ is detected, the economic system will now use the model $Y_t=X'_t\beta_1 + \xi_t, ~ t=k+1,\cdots $
    and monitoring will continue to ensure that this new model can explain the new data that will arrive.
  \end{exam}
  The sequential change-point problem presented in Example \ref{ex1} is often designed by the terminology "fault detection control" or "engineering process control".
  The aim of the control scheme is to trigged an alarm and stop the system when disorder is detected. Since the diagnostic, the resetting and the malfunctioning can
  have a cost, the scheme will be able to reduce as soon as possible the false alarm and the delay of detecting disorder.
  For more review of this approach, see Pollak and Siegmund (1991), Basseville and Nikiforov (1993),  Lai (1998), Lai and Shan (1999), Mei (2006),
   Moustakides (2008) $\cdots$\\
   The problem presented in Example \ref{ex2} is designed by the terminology "monitoring scheme" (see for instance Chu {\it et al.} \cite{Chu1996}) ; the term "on-line segmentation" or "sequential segmentation"
   is often used in application in signal processing (see for instance Basseville and Nikiforov \cite{Basseville1993} p. 2). These terminologies refer to the
   sequential procedure where false alarms are less crucial (than in the previous example, because here, the system will deal with in the next stage); in such case, it is important to estimate the model described by the system after change.\\
    This paper focus on such sequential change-point problem and we will follow the paradigm of Chu {\it et al.} \cite{Chu1996},
   Berkes {\it et al.}  (2004), Horv\'ath {\it et al.} (2004), Zeileis {\it et al.} (2005),
   Aue {\it et al.} (2009), Na {\it et al.}  (2011) $\cdots$, who have seen this problem as a classical hypothesis testing with a fixed probability of type I error.\\
  An important turning on this topic was made with the paper of Chu {\it et al.}. They considered the sequential change in regression model and pointed out the effects of repeating retrospective
   test when new data are observed; this can increase the probability of type I error of the test.
   They successfully applied a fluctuation test to solve the sequential change-point problem.
  Two procedures are developed based on cumulative sum (CUSUM) of residuals and recursive parameter fluctuations.
   Their idea has been generalized and several procedures are now based on this approach.
    Leisch {\it et al.} \cite{Leisch2000} introduced the generalized fluctuation test based on the recursive moving estimator which contains the test of Chu {\it et al.} \cite{Chu1996} as a special case.
    Horváth {\it et al.} \cite{Horvath2004} introduced residual CUSUM monitoring procedure where the recursive parameter is based on the historical data.
    This procedure has been generalized by Aue {\it et al.} \cite{Aue2006} to the class of linear model with dependent errors.
     Berkes {\it et al.} \cite{Berkes2004} considered sequential changes in the parameters of GARCH process.
    According to the fact that the functional limit theorem assumed by  Chu {\it et al.} \cite{Chu1996} is not satisfied by the squares of residuals of GARCH process, they developed a procedure based on quasi-likelihood scores.
     Na {\it et al.} \cite{Lee2011} developed a monitoring procedure for the detection of parameter changes in general time series models. They show that under the null hypothesis of no change, their detector statistic converges weakly to a known distribution. However, the asymptotic behavior of this detector is   unknown under the alternative of parameter changes.\\

 In this new contribution, we consider a large class of causal time series and investigate  the asymptotic behavior under the null hypothesis of no change but also under the alternative hypothesis of change. More precisely, let $M , f : \R^{\N} \rightarrow \R$ be
 measurable functions,  $(\xi_t)_{t\in \Z}$ be a sequence of
 centered independent and identically distributed (iid)  random variables satisfying $\textrm{var}(\xi_0)= \sigma^2$  and let $\Theta $ be a fixed compact subset of $\R^{d}$. Let $T \subset \Z$, and for any $ \theta \in \Theta$, define ~\\

 {\bf Class $\mathcal{M}_T(M_{\theta},f_{\theta})$:} {\it The process
 $X=(X_t)_{t\in\Z}$ belongs to $\mathcal{M}_T(M_{\theta},f_{\theta})$ if it satisfies the relation:}
 \begin{equation}\label{model}
 X_{t+1}=M_{\theta} \big ((X_{t-i})_{i\in \N}\big ) \xi_t+f_{\theta}\big
 ((X_{t-i})_{i\in \N}\big )\quad\mbox{for all $t\in T$}.
 \end{equation}

    \noindent The existence and properties of this general class of causal and affine processes were studied in Bardet and Wintenberger \cite{Bardet2009}. Numerous classical time series (such as AR$(\infty)$, ARCH$(\infty)$, TARCH$(\infty)$, ARMA-GARCH or bilinear processes) are included in ${\cal M}_\Z(M,f)$. The off-line change detection for such class of models has already been studied in Bardet {\it et al.} \cite{Bardet2010}  and Kengne \cite{Kengne2011}.\\

Suppose now that we have observed $X_1,\cdots,X_n$ which are available historical data such that there exists $ \theta^*_0 \in \Theta$ such as $(X_1,\cdots,X_n)$ belongs to  $\mathcal{M}_{\{1,\cdots,n\} }(M_{\theta^*_0},f_{\theta^*_0})$.
  Then, we observe new data $ X_{n+1},X_{n+2}\cdots,X_k,\cdots $: the monitoring scheme starts.
  For each new observation, we would like to know if a change occurs in the parameter $\theta^*_0$.
  More precisely, we consider the following test :\\

 \noindent  $\mathbf{H_0}$: $\theta^*_0$ is constant over the observation $X_1,\cdots,X_n,X_{n+1},\cdots$ {\it i.e.} the observations  $X_1,\cdots,X_n,X_{n+1},\cdots$ belong to
           $\mathcal{M}_{\N }(M_{\theta^{*}_0},f_{\theta^*_0})$;\\
 ~\\
 $\mathbf{H_1}$ : there exist $k^*>n$, $\theta^*_1 \in \Theta$ such that $X_1,\cdots,X_n,X_{n+1},\cdots,X_{k^*},X_{k^*+1},\cdots$ belongs to\\
 $\mathcal{M}_{\{1,\cdots,k^*\} }(M_{\theta^*_0},f_{\theta^*_0})\bigcap \mathcal{M}_{\{k^*+1,\cdots\} }(M_{\theta^*_1},f_{\theta^*_1})$.\\

  When new data arrive, Chu {\it et al.} \cite{Chu1996} proposed to compute an estimator of the parameter based on all the observations and to compare it to
  an estimator based on historical data. A large distance between both these estimators means that new data come from a model with different parameters.
  Then the  null hypothesis $H_0$ is rejected and the new parameter is considered; otherwise, the monitoring scheme continues.
  In their procedure, Leisch {\it et al.} \cite{Leisch2000} suggested to compute the recursive estimators on a moving window with a fixed width.
  They fixed a monitoring horizon so that, the procedure will stop after a fixed number of steps even if no change is detected.
  As Chu {\it et al.} \cite{Chu1996}, the recursive estimators computed by Na {\it et al.} \cite{Lee2011} are based on all the observations.
  As we will see in the next sections, their procedure cannot be effective in terms of detection delay or if a small change in the parameter occurs.
  \\
  \indent For any $k \geq 1$,  $\ell,\ell'\in \{1,\cdots,k \}$ (with $\ell\leq \ell'$) let $\widehat{\theta}(X_{\ell},\cdots,X_{\ell'})$ be the quasi-maximum likelihood estimator (QMLE in the sequel) of the parameter computed on $\{\ell,\cdots,\ell' \}$ as it is  defined in \eqref{theta}.
  When new data arrive at time $k\geq n$, we explore the segment $ \{ \ell, \ell+1,\cdots,k \} $ with  $\ell\in \{n-v_n, n-v_n+1 \cdots k-v_n \}$ (where $(v_n)_{n\in \N}$ is a fixed sequence of integer numbers) that the distance between $\widehat{\theta}(X_{\ell},\cdots,X_{k})$ and $\widehat{\theta}(X_1,\cdots,X_{n})$  is the largest. If the norm $ \| \widehat{\theta}(X_{\ell},\cdots,X_{k}) - \widehat{\theta}(X_{1},\cdots,X_{n})\|$ is greater than a suitable critical value, then
  $H_0$ is rejected and a model with a new parameter is considered; otherwise, the monitoring scheme continues.
  More precisely, we construct a detector that takes into account the distance between  $\widehat{\theta}(X_{\ell},\cdots,X_{k})$ and $\widehat{\theta}(X_1,\cdots,X_{n})$.
  It is shown that this detector is almost surely finite under the null hypothesis and almost surely diverges to infinity under the alternative.
  Hence, the consistency of our procedure follows.\\
  Finally, Monte-Carlo experiments have been done, comparing our procedure to the ones of Horváth {\it et al.} \cite{Horvath2004} (see also Aue {\it et al.} \cite{Aue2006}) and Na {\it et al.} \cite{Lee2011}. It appears that our procedure outperforms these other procedures in terms of test power and detection delay in different frames. An application to financial data (Nikkei 225, S$\&$P 500 and FTSE 100 stock index) allows to detect changes in these data in accordance with
  historical and economic events. \\

  In the forthcoming Section 2 the assumptions and the definition of the quasi-likelihood estimator are provided. In Section 3 we present the monitoring procedure and the asymptotic results.
  Section 4 is devoted to a simulation study for  AR$(1)$ and GARCH$(1,1)$ processes. In Section 5 we apply our procedure to famous financial data. The proofs of the main results are provided in Section 6.

\section{Assumptions and definition of the quasi-likelihood estimator}

\subsection{Assumptions on the class of models ${\cal M}_\Z (f_\theta,M_\theta)$}

Let $\theta\in\R^d$ and $M_\theta$ and $f_\theta$ be numerical functions such that for all $(x_i)_{i\in \N}\in \R^{\N}$,
 $M_\theta\big((x_i)_{i\in \N}\big)\neq 0$ and $f_\theta\big((x_i)_{i\in \N}\big) \in \R$. Denote $ h_\theta := M_\theta^2$.  We will use the following
 classical notations:
 \begin{enumerate}
 \item  $\|\cdot\|$ applied to a vector denotes the Euclidean norm of
    the vector;
    \item for any compact set $ \mathcal{K} \subseteq\R^d$ and for any
    $g:\mathcal{K} \longrightarrow\R^{d'}$, $ \|g\|_\mathcal{K} =\sup_{\theta\in \mathcal{K} }(\|g(\theta)\|)$;
    \item for any set $ \mathcal{K} \subseteq\R^d$,  $\overset{\circ}{\mathcal{K} }$ denotes the interior of $\mathcal{K} $.
 \end{enumerate}

 \noindent Throughout the sequel, we will assume that the functions $ \theta \mapsto M_\theta$ and $\theta \mapsto f_\theta$
 are twice continuously differentiable on $\Theta$.
 Let  $\Psi_\theta=f_\theta, \, M_\theta$ and $i=0,\, 1, \, 2$, then for any compact set $\mathcal{K} \subset \Theta$ define \\
~\\
{\bf Assumption A$_i(\Psi_\theta,\mathcal{K} )$}: {\em Assume that
$\|{\partial^i\Psi_\theta(0)}/{\partial\theta^i}\|_\Theta<\infty$
 and there exists a sequence of non-negative real numbers $(\alpha^{(i)}_j(\Psi_\theta,\mathcal{K} ))_{j\geq 1}$ such that $ \sum\limits_{j=1}^{\infty} \alpha^{(i)}_j(\Psi_\theta,\mathcal{j} ) <\infty$ and
\begin{equation*}
\Big\|\dfrac{\partial^i\Psi_\theta(x)}{\partial\theta^i}-\dfrac{\partial^i\Psi_\theta(y)}{\partial\theta^i}\Big\|_\mathcal{K}
\leq
\sum\limits_{j=1}^{\infty}\alpha^{(i)}_j(\Psi_\theta,\mathcal{K}
)|x_j-y_j|\quad \mbox{for all}~x, y \in \R^{\N}.
\end{equation*}}
\noindent In the sequel we refer to the particular case called "ARCH-type process"
 if $f_\theta =0$ and if the following assumption holds with $ h_\theta = M_\theta^2$:\\
~\\
{\bf Assumption A$_i(h_\theta,\mathcal{K} )$}: {\em Assume that
$\|{\partial^i h_\theta(0)}/{\partial\theta^i}\|_\Theta<\infty$
 and there exists a sequence of non-negative real numbers $(\alpha^{(i)}_j(h_\theta,\mathcal{K} ))_{j\geq 1}$ such as $ \sum\limits_{j=1}^{\infty} \alpha^{(i)}_j(h_\theta,\mathcal{K} ) <\infty$ and
\begin{equation*}
\Big\|\dfrac{\partial^i
h_\theta(x)}{\partial\theta^i}-\dfrac{\partial^i
h_\theta(y)}{\partial\theta^i}\Big\|_\mathcal{K}
\leq \sum\limits_{j=1}^{\infty}\alpha^{(i)}_j(h_\theta,\mathcal{K} )|x^2_j-y^2_j|\quad \mbox{for all}~x, y \in \R^{\N}.\\
\end{equation*}}
The Lipschitz-type hypothesis $A_i(\Psi_\theta,\mathcal{K} )$ are classical when studying the existence of solutions of the general model (see for instance \cite{Doukhand}). Using a result of \cite{Bardet2009}, for each model $\mathcal{M}_{\Z }(M_{\theta},f_{\theta})$ it is interesting to define the following set:
   \begin{multline*}
  \Theta(r):=\Big \{\theta \in \Theta,\, A_0(f_\theta,\{\theta\})\, \,\mbox{and} \,A_0(M_\theta,\{\theta\}) \, \textrm{hold with}\,
 \sum\limits_{j \geq 1}^{ } \alpha^{(0)}_j(f_\theta,\{\theta\}) + (\E| \xi_0 |^r)^{1/r}\sum\limits_{j \geq 1}^{ } \alpha^{(0)}_j(M_\theta,\{\theta\}) <1 \Big \}~~~~ \\
  \bigcup \Big \{\theta \in \Theta,\, f_\theta = 0  \textrm{ and}\,  \,A_0(h_\theta,\{\theta\}) \, \textrm{holds with}\,
  (\E| \xi_0 |^r)^{2/r}\sum\limits_{j \geq 1}^{ } \alpha^{(0)}_j(h_\theta,\{\theta\}) <1 \Big \}.
 \end{multline*}
Then, if $\theta\in\Theta(r)$ the existence of a unique causal, stationary and ergodic  solution  $X=(X_t)_{t\in \Z}\in {\cal
 M}_\Z(f_\theta,M_\theta)$ is ensured (see more details in \cite{Bardet2009}). The subset $\Theta(r)$ is defined as a union to consider
 accurately general causal models and ARCH-type models simultaneously.\\
~\\
 Here there are assumptions required for studying QMLE asymptotic properties:\\
 ~\\
 \noindent {\bf Assumption D$(\Theta)$:} $\exists\underline{h}>0$ such that
$\displaystyle \inf_{\theta \in
 \Theta}(|h_\theta(x)|)\geq \underline{h}$ for all $x\in \R^{\N}.$ \\
~\\
 {\bf Assumption Id($\Theta$):} For all  $(\theta,\theta')\in \Theta^2$,
$$ \Big( f_{\theta}(X_0,X_{-1},\cdots)=f_{\theta'}(X_0,X_{-1},\cdots)~\mbox{and}~h_{\theta}(X_0,X_{-1},\cdots)=h_{\theta'}(X_0,X_{-1},\cdots) \ \text{a.s.}\Big) \Rightarrow \ \theta = \theta'.$$
{\bf Assumption Var($\Theta$):} For all $\theta  \in \Theta $, one
 of the families $ \big( \dfrac{\partial f_{\theta}}{\partial \theta^{i}}(X_0,X_{-1},\cdots) \big)_{1\leq i \leq d} \quad
  \mbox{or} \quad \big( \dfrac{\partial h_{\theta}}{\partial \theta^{i}}(X_0,X_{-1},\cdots) \big)_{1\leq i \leq d}  \quad $
 is a.s. linearly independent.\\
~\\
\textbf{Assumption K}($f_\theta,M_\theta,\Theta$): for i= 0, 1, 2, \textbf{A}$_i(f_\theta,\Theta)$ and \textbf{A}$_i(M_\theta,\Theta)$
 (or \textbf{A}$_i(h_\theta,\Theta)$) hold and there exists $\ell>2$ such that $ \alpha_j^{(i)}(f_{\theta},\Theta) + \alpha_j^{(i)}(M_{\theta},\Theta) + \alpha_j^{(i)}(h_{\theta},\Theta) =
      \mathcal{O}(j^{-\ell})$ for $j\in \N$.\\
     ~\\
Note that in this last assumption, as in \cite{Bardet2009}, we use the convention that if \textbf{A}$_i(M_\theta,\Theta)$ holds then $\alpha_\ell^{(i)}( h_\theta,\Theta) = 0$ and if
     \textbf{A}$_i(h_\theta,\Theta)$ holds then $\alpha_\ell^{(i)}( M_\theta,\Theta) = 0$.
 \noindent

 \subsection{Two first examples}

 \begin{enumerate}
    \item ARMA$(p,q)$ processes.\\
          Consider the ARMA$(p,q)$ process defined by:
          \begin{equation}\label{arma}
           X_t+\sum_{i=1}^p a^*_i X_{t-i}  = \sum_{j=0}^q b_j^* \xi_{t-j}  ~ , ~ t \in \Z
           \end{equation}
 with $b_0^*\neq 0$, $ \theta^*_0=(a_1^*,\cdots,a_p^*,b_0^*,\cdots,b_q^*) \in \Theta \subset \R^{p+q+1}$ and $(\xi_t)$ a white noise such as $\E(\xi_0^2)=1$. When $\sum_{j=0}^q b_j^* X^j \neq 0$ and $1+\sum_{i=0}^p a_i^* X^i \neq 0$ for all $|X|\leq 1$, this process can be also written as:
 $$
 X_t=b_0^* \xi_t + \sum_{j=1}^\infty \phi_j(\theta_0^*) \,  X_{t-i}~ , ~ t \in \Z
 $$
 where $\theta\in \Theta \mapsto \phi_j(\theta)$ are functions only depending on $\theta$ and decreasing exponentially fast to $0$  ($j\to \infty$).
  The process \eqref{arma} belongs to the class ${\cal M}_{\Z}(M_{\theta^*_0},f_{\theta^*_0})$ where
  $f_\theta(x_1,\cdots)= \sum_{j\geq 1} \phi_j(\theta)x_j $ and $ M_\theta \equiv b_0^* ~ \text{for all}~ \theta \in \Theta$.
  Then Assumptions {\bf D}($\Theta$), {\bf A}$_0(f_\theta,\Theta)$, {\bf A}$_0(M_\theta,\Theta)$ hold  with $\underline{h}=|b_0^*|>0$ and $\alpha^{(0)}_j(f_\theta,\Theta) = \|\phi_j(\theta)\|_\Theta$ while $\alpha^{(0)}_j(M_\theta,\Theta) =0$ for $j\in \N^*$.
   Assumption \textbf{K}($f_\theta,M_\theta,\Theta$) holds since there exists $c>0$ and $C>0$ such as $|\phi_j|\leq C \, e^{-cj}$ for $j\in \N$. Moreover, if $(\xi_t)$ is a sequence of non-degenerate random variables ({\it i.e.} $\xi_t$ is not equal to a constant), Assumptions
   {\bf Id}($\Theta$) and {\bf Var}($\Theta$) hold. Finally, for any $r\geq 1$ such that $ \E|\xi_0|^r<\infty$, then
$$\Theta(r)=\Big \{\theta \in \R^{p+q+1},~ \sum_{j\geq 1} |\phi_j(\theta)|  <1 \Big  \}.$$
Note that if $\theta \in \Theta(r)$ with $r\geq 1$ then the previous conditions of stationarity $\sum_{j=0}^q b_j X^j \neq 0$ and $1+\sum_{i=0}^p a_i X^i \neq 0$ for all $|X|\leq 1$ are satisfied.

    \item GARCH$(p,q)$ processes.\\
          Consider the GARCH$(p,q)$ process defined by:
 \begin{equation}  \label{garch}
 X_t= \sigma_t \, \xi_t \ ,\ \sigma_t^2 = \alpha_{0}^*+ \sum^{p}_{j=1} \alpha_{j}^*X^2_{t-j} + \sum^{q}_{j=1} \beta_{j}^*\sigma^2_{t-j}  \ ,  \ \  t\in \Z
 \end{equation}
 with $\E(\xi_0^2)=1$ and $ \theta^*_0:= (\alpha^*_0,\cdots,\alpha^*_p,\beta^*_1,\cdots,\beta^*_q ) \in \Theta $ where $\Theta$ is a
 compact subset of $ ]0 , \infty[\times [0 , \infty [^{p+q}$   such that
 $\sum_{j=1}^{p}\alpha_j +  \sum_{j=1}^{q}\beta_j < 1 $ for all $\theta \in \Theta$.
 %
 Then there exists (see Bollerslev \cite{Bollerslev1986} or Nelson and Cao \cite{Nelson1992}) a nonnegative sequence $(\psi_j(\theta^*_0))_{j\geq 0}$
 such that $\sigma_t^2 = \psi_0(\theta^*_0) + \sum^{}_{j\geq 1} \psi_j(\theta^*_0)X^2_{t-j}$  with  $ \psi_0(\theta^*_0) = \alpha^*_0/(1-\sum^{q}_{j=1} \beta_{j}^*)$.\\
 This process belongs to the class ${\cal M}_{\Z}(M_{\theta^*_0},f_{\theta^*_0})$ where $ f_\theta \equiv 0$ and
  $M_\theta(x_1,\cdots)= \sqrt{\psi_0(\theta) + \sum_{j\geq 1} \psi_j(\theta)x^2_j} $  for all $\theta \in \Theta$.
   Assumption {\bf D}($\Theta$)  holds with $\underline{h}= \underset{\theta\in \Theta} {\mbox{inf}}(\psi_0(\theta))>0$.
  If there exists $0<\rho_0<1$ such that for any $\theta \in \Theta, ~ \sum_{j=1}^{q}\alpha_j +  \sum_{j=1}^{p}\beta_j \leq \rho_0$ then
 the sequences $( \|\psi_j(\theta)\|_\Theta)_{j\geq 1}$, $( \|\psi'_j(\theta)\|_\Theta)_{j\geq 1}$ and $( \|\psi''_j(\theta)\|_\Theta)_{j\geq 1}$
  decay exponentially fast (see Berkes {\it et al.} \cite{Berkes2003}) and Assumption \textbf{K}($f_\theta,M_\theta,\Theta$) holds.
  Moreover, $(\xi^2_t)$ is a sequence of non-degenerate random variables ({\it i.e.} $\xi^2_t$ is not equal to a constant), Assumptions
   {\bf Id}($\Theta$) and {\bf Var}($\Theta$) hold.  Finally for $r\geq 2$ we obtain
   $$  \Theta(r)= \Big\{  \theta \in \Theta ~ ; ~  (\E|\xi_0|^r)^{2/r} \sum_{j=1}^{\infty }\phi_j(\theta) < 1 \Big\}  .$$

 \end{enumerate}

  \subsection{The quasi-maximum likelihood estimator}

  Let $k \geq n\geq 2$, if $ (X_1, \cdots, X_k) \in \mathcal{M}_{ \{1,\cdots,k\} }(M_\theta,f_\theta)$, then
        for $ T\subset \{1,\cdots,k \} $,  the conditional quasi-(log)likelihood computed on $T$ is given by:
\begin{equation}\label{defML}
 L(T,\theta):= -\dfrac{1}{2} \sum\limits_{t \in T} q_t(\theta)   ~~ \text{with}  ~~ q_t(\theta)= \dfrac{(X_t-f^t_{\theta})^2}{h_\theta^t} + \log(h_\theta^t)
 \end{equation}
       where $ f^t_{\theta}=f_{\theta}\big(X_{t-1},X_{t-2}\ldots\big) $, $ M^t_{\theta}=M_{\theta}\big(X_{t-1},X_{t-2}\ldots\big) $ and
   $ h^t_{\theta}= {M^{t}_{\theta}}^2 $. The classical approximation of this conditional log-likelihood (see more details in Bardet and Wintenberger \cite{Bardet2009}) is given by:

   \begin{equation}\label{hatML}
   \widehat L(T,\theta):=-\frac{1}{2}\sum\limits_{t \in T}\widehat q_t(\theta)\quad
 \textrm{where}\;\;\;\widehat{q}_t(\theta):=\frac{\big(X_t-\widehat{f}^t_{\theta}\big)^2}{\widehat{h}{^t_{\theta}}} +\log\big(\widehat{h}{^t_{\theta}}\big)
 \end{equation}
 with
 $ \widehat{f}^t_{\theta}=f_{\theta}\big(X_{t-1},\ldots,X_{1},0,0,\cdots\big)$, $ \widehat{M}^t_{\theta}=
 M_{\theta}\big(X_{t-1},\ldots,X_{1},0,0,\cdots\big) $ and $\widehat{h}^t_{\theta}=( \widehat{M}^{t}_{\theta} )^2$.\\
   For $ T\subset \{1,\cdots,k \} $, define the quasi maximum-likelihood estimator (QMLE) of $\theta$ computed on $T$ by
   \begin{equation}\label{theta}
   \widehat{\theta}(T):=  \underset{\theta\in \Theta} {\mbox{argmax}}(\widehat{L}(T,\theta)).
   \end{equation}
  In Bardet and Wintenberger \cite{Bardet2009} it was established that if $(X_1,\cdots,X_n)$ is an observed trajectory of $X \subset {\cal M}_{\Z}(f_{\theta_0^*},M_{\theta_0^*})$ with  $\theta^*_0 \in \overset{\circ}{\Theta}(4)$ and if $\Theta$ is a compact set such as Assumptions \textbf{A}$_i$($f_\theta,M_\theta,\Theta$) (or \textbf{A}$_i$($h_\theta,\Theta$)) hold for $i=0,1,2$ and under Assumptions \textbf{D}$(\Theta)$, \textbf{Id}($\Theta$),
 \textbf{Var}($\Theta$), \textbf{K}($f_\theta,M_\theta,\Theta$), then
 \begin{equation}\label{TLC}
 \sqrt n \big (\widehat{\theta}(T_{1,n})-\theta_0^* \big) \limiteloin {\cal N}\big ( 0\, , \,F\,{G}^{-1}\,F \big ),
 \end{equation}
with
\begin{equation}\label{FG}
 \displaystyle G:= \E \Big [ \dfrac{\partial q_0(\theta^*_0)}{\partial \theta} \dfrac{\partial q_0(\theta^*_0)}{\partial \theta} ' \Big] \quad \mbox {and}\quad \displaystyle  F:=  \E \Big [ \dfrac{\partial^2 q_0(\theta^*_0)}{\partial \theta \partial \theta'}  \Big],
  \end{equation}
 where $'$ denotes the transpose and with $q_0$ defined in \eqref{defML}. Note that under assumptions \textbf{D}$(\Theta)$ and \textbf{Var}$(\Theta)$, $G$ is symmetric positive definite (see \cite{Kengne2011}) and $F$ is non-singular (see \cite{Bardet2009}). Also define the matrix
 \begin{equation}\label{hatFG}
 \widehat{G}(T):= \dfrac{1}{\mbox{Card}(T)} \sum \limits_{t \in T} \Big( \dfrac{\partial \widehat{q}_{t}( \widehat{\theta}(T))}{\partial \theta} \Big)
                       \Big( \dfrac{\partial \widehat{q}_{t}( \widehat{\theta}(T))}{\partial \theta} \Big)'\quad\mbox{and}\quad
    \widehat{F}(T):= - \dfrac{2}{\mbox{Card}(T)} \Big(  \dfrac{\partial^2 \widehat{L}_{m} (T,\widehat{\theta}(T))}{\partial \theta \partial \theta'}
   \Big).
   \end{equation}
    Under the previous assumptions, $\widehat{G}(T_{1,n})$ and $\widehat{F}(T_{1,n})$  converge almost surely to $G$ and $F$ respectively.
    Hence,  \begin{equation}\label{central}
\sqrt n \,  \widehat{G}(T_{1,n})^{-1/2} \widehat{F}(T_{1,n})\big (\widehat{\theta}(T_{1,n})-\theta_0^* \big) \limiteloin {\cal N}\big ( 0\, , \,I_d \big )
            \end{equation}
 with $I_d$ the identity matrix.  This result will be the starting point of the following monitoring procedure.

 \section{The monitoring procedure and asymptotic results}

   \subsection{The monitoring procedure}\label{monit_proc}
    In the sequel, $(X_1,\cdots,X_n)$ is supposed to be the historical available observations belonging to the class ${\cal M}_{\{1,\cdots,n\}}(f_{\theta_0^*},M_{\theta_0^*})$.
     For $ 1 \leq \ell \leq \ell' $, denote
     $$
     T_{\ell,\ell'}:= \{\ell,\ell+1,\cdots,\ell' \}.
     $$
   At a monitoring instant $k$, our procedure evaluates the difference between $\widehat{\theta}(T_{\ell,k})$ and  $\widehat{\theta}(T_{1,n})$ for any $\ell=n,\cdots,k$.
   More precisely, from \eqref{central}, for any $k>n$ define the statistic (called the detector)
    $$ \widehat{C}_{k,\ell}:=  \sqrt{n}\, \dfrac{k-\ell}{k}\, \big\|\widehat{G}(T_{1,n})^{-1/2} \widehat{F}(T_{1,n})\big(\widehat{\theta}(T_{\ell,k}) - \widehat{\theta}(T_{1,n})\big)\big\|$$
    for $\ell=n,\cdots,k$. Since the matrix $\widehat{G}(T_{1,n})$ is  asymptotically symmetric and positive definite (see \cite{Kengne2011}),
    $ \widehat{G}(T_{1,n})^{-1/2}$ exists for $n$ large enough and $\widehat{C}_{k,\ell}$ is well defined.
    At the beginning of the monitoring scheme and   when $\ell$  is close to $k$, the length of $T_{\ell,k}$ is too small, therefore
    the numerical algorithm used to compute $\widehat{\theta}(T_{\ell,k})$ cannot converge. This can introduce a large distortion in the procedure.
    To avoid this, we introduce a  sequence of integer numbers $(v_n)_{n \in \N}$ with $v_n<< n$ and compute
    $\widehat{C}_{k,\ell}$ for $\ell \in \{ n-v_n, n-v_n+1,\cdots,k-v_n\}$. Thus, for any $k>n$ denote
    $$  \Pi_{n,k} :=   \{ n-v_n, n-v_n+1,\cdots,k-v_n\} .$$
    For technical reasons, assume that,
    $$ v_n \to \infty \quad \mbox{and}\quad v_n/{\sqrt{n}} \to 0 ~~(n\to \infty).$$

 \noindent According to Remark 1 of \cite{Kengne2011}, we can choose $ v_n= [(\log n)^{\delta}]$ with $\delta >1$.\\
    Note that, if change does not occur at time $k>n$, for any $\ell \in  \Pi_{n,k}$, the two estimators $\widehat{\theta}(T_{\ell,k})$ and $\widehat{\theta}(T_{1,n})$
    are close and the detector $\widehat{C}_{k,\ell}$ is not too large.
    Hence, the monitoring scheme rejects $H_0$ at the first time $k>n$ where there exists $\ell \in \Pi_{n,k}$ satisfying  $\widehat{C}_{k,\ell}>c$ for a fixed constant $c>0$.
  To be more general, we will use a $b:(0,\infty)\mapsto (0,\infty)$, called a boundary function satisfying:\\
  ~\\
 {\bf Assumption B:}  $b:(0,\infty)\mapsto (0,\infty)$ is a
  non-increasing and continuous function such as $ \underset{0<t<\infty} {\mbox{Inf}}b(t)>0$.\\
  ~\\
  Then the monitoring scheme rejects $H_0$ at the first time $k>n$ such as there exists $\ell \in \Pi_{n,k}$ satisfying  $\widehat{C}_{k,\ell}>b((k-\ell)/n)$.
  Hence define the stopping time:
  $$\tau(n):= \text{Inf}\Big  \{ k>n ~ / ~ \exists \ell \in \Pi_{n,k}, ~ \widehat{C}_{k,\ell}>b((k-\ell)/n)\Big  \}
  = \text{Inf} \Big\{ k>n ~ / ~   \underset{\ell \in \Pi_{n,k}} {\mbox{max}} \dfrac{\widehat{C}_{k,\ell}}{b((k-\ell)/n)}>1 \Big\} .$$
 Therefore, we have
 \begin{equation}\label{prob_stop_tim}
 P\{ \tau(n)<\infty \}=P\Big\{ ~ \underset{\ell \in \Pi_{n,k}} {\mbox{max}} \dfrac{\widehat{C}_{k,\ell}}{b((k-\ell)/n)}>1 ~ \text{for some} ~ k>n \Big \}
   = P\Big\{ ~ \underset{k>n} {\mbox{sup}} ~ \underset{\ell \in \Pi_{n,k}} {\mbox{max}} \dfrac{\widehat{C}_{k,\ell}}{b((k-\ell)/n)}>1 \Big \} .
  \end{equation}
  The challenge is to choose a suitable boundary function $b(\cdot)$ such as for some given $\alpha \in (0,1)$
  $$ \lim_{n \rightarrow \infty} P_{H_0}\{ \tau(n)<\infty \}=\alpha$$
  and
  $$ \lim_{n \rightarrow \infty}  P_{H_1}\{ \tau(n)<\infty \}=1$$
 where the hypothesis $H_0$  and $H_1$ are specified in Section \ref{intro}. \\In the case where $b(\cdot)$ is a constant positive value,
 $b \equiv c$ with $c>0$, these conditions lead to compute a threshold $c=c_\alpha$ depending on $\alpha$.
 If change is detected under $H_1$ i.e.  $ \tau(n)<\infty $ and $\tau(n)>k^*$, then the detection delay is defined by
 $$\widehat{d}_n= \tau(n) - k^*.$$
   Using the previous  notations,  Na {\it et al.} \cite{Lee2011} used the following detector
    $$
    \widehat{D}_k:=\sqrt{n} \big\|\widehat{G}(T_{1,n})^{-1/2} \widehat{F}(T_{1,n})\big(\widehat{\theta}(T_{1,k}) - \widehat{\theta}(T_{1,n})\big)\big\| .
    $$
   At the step $k$ of the monitoring scheme, their recursive estimator is based on
   $X_1,\cdots,X_n,\cdots,X_{k}$. One can see that this estimator is highly influenced by the historical data.
   Assume that a change occurs at time $k^*\leq k$, in the sequel of the procedure, the recursive estimator contents the observations
    $X_1,\cdots,X_n,\cdots,X_{k^*-1}$ which depends on $\theta^*_0$.
   Then, one must wait longer before the difference between $\widehat{\theta}(X_1,\cdots,X_n)$ and $\widehat{\theta}(X_1,\cdots,X_n,\cdots,X_{k})$
   becomes significant at a step $k>k^*$. Therefore, their procedure cannot be effective in terms of detection delay.
   Moreover, if $n$ tends to infinity, it is not almost sure that this change will be detected.
   These are confirmed by the results of simulations (see Section 4).\\
   \indent  Berkes {\it et al.} (2004) considered an estimator based on historical data to compute
   the quasi-likelihood scores.  They used the fact that the partial derivatives applied to a vector $\textbf{u}$ is equal to $0$ if and only if $\textbf{u}$
   is the true parameter of the model.   So, when change occurs, their detector growths asymptotically to infinity.
   Therefore, their procedure is consistent. They proved this result for GARCH(p,q) models.

   \subsection{Asymptotic behaviour under the null hypothesis}
  Under $H_0$, the parameter $\theta^*_0$ does not change over the new observations.
  Thus we have the result
\begin{Theo}\label{theo1}
 Assume \textbf{D}$(\Theta)$, \textbf{Id}($\Theta$),
 \textbf{Var}($\Theta$), \textbf{K}($f_\theta,M_\theta,\Theta$), {\bf B} and $ \theta^*_0 \in \overset{\circ}{\Theta}(4)$. Under the null hypothesis $H_0$, then
  $$
  \lim_{n \rightarrow \infty} P\{\tau(n)<\infty \} = P\Big\{\sup_ {t>1} \sup_{1<s<t} \dfrac{\| W_d(s) - s W_d (1))\|}{t~b(s)}>1 \Big\}
  $$
where $W_d$ is a $d$-dimensional standard Brownian motion.
\end{Theo}
 In the simulations, we will use the most ``natural'' boundary function $b(\cdot)=c$ with $c$ a positive constant since it satisfies the above assumptions imposed to $b(\cdot)$. In such case, the forthcoming corollary indicates that the asymptotic distribution of Theorem \ref{theo1} can be easily computed:
 \begin{cor}\label{cor1}
  Assume $b(t)=c>0$ for $t\geq 0$. Under the assumptions of Theorem \ref{theo1},
  $$
  \lim_{n \rightarrow \infty} P\{\tau(n)<\infty \}  = P\Big\{\sup_ {t>1} \sup_{1<s<t} \frac{1}{t} \| W_d(s) - s W_d (1))\| >c \Big\}
    = P\{  U_d > c \} $$
  where
  $U_d= \sup_{0<u<1 } ~ f(u) \,  \| W_d(u)\|
 \quad \mbox{with}\quad f(u)=\dfrac {\sqrt{9-u}+\sqrt{1-u}}{\sqrt{9-u}+3\sqrt{1-u}}\, \Big ( \dfrac { 2}  {3-u+\sqrt{(9-u)(1-u)}}\Big )^{1/2}.$

 \end{cor}

  \begin{rem}
   ~
   \begin{enumerate}
  \item By the law of the iterated logarithm, it comes that
  $$ \sup_{1<s<t} \frac{1}{t} \| W_d(s) - s W_d (1))\| ~ \overset{ \text{a.s.}}{\underset{t \rightarrow \infty} \longrightarrow} ~ \|W_d (1))\| .$$
  So, the two distributions $\sup_{1<s<t} \frac{1}{t} \| W_d(s) - s W_d (1))\|$ as $t\rightarrow \infty$ (resp. $ t \rightarrow 1$) and
  $f(u) \,  \| W_d(u)\|$ as $u\rightarrow 1$ (resp.  $u \rightarrow 0$) are equal.
  It is easy to show (see proof of Corollary \ref{cor1}) that
   $$\sup_ {t>1} \sup_{1<s<t} \frac{1}{t} \| W_d(s) - s W_d (1))\|  ~ \overset{ \mathcal{D}}{  =} ~  \sup_{0<u<1 } ~ f(u) \,  \| W_d(u)\|.$$
  \item  Under the null hypothesis, it holds that $ \widehat{\theta}(T_{1,n}) \limitepsn \theta^*_0 $ (see \cite{Bardet2009}). Thus denote
  $$ \widehat{C}^{(0)}_{k,\ell}:=  \sqrt{n}\, \dfrac{k-\ell}{k}\, \big\|\widehat{G}(T_{1,n})^{-1/2} \widehat{F}(T_{1,n})\big(\widehat{\theta}(T_{\ell,k}) -  \theta^*_0 \big)\big\|.$$
  Under the assumptions of Theorem \ref{theo1}, one can easily show that
 $$ \sup_{k>n}  ~\max_{ \ell \in \Pi_{n,k}} ~ \dfrac{1}{b((k-\ell)/n)}\, \big| \widehat{C}_{k,\ell} - \widehat{C}^{(0)}_{k,\ell} \big|=o_P(1) ~ ~  \text{as} ~ n\rightarrow \infty.$$
 Thus, the Theorem \ref{theo1} still holds if the stopping time $\tau(n)$ is computed by using the detector $\widehat{C}^{(0)}_{k,\ell}$.
 Hence, if the parameter $\theta^*_0$ of the historical observations is known, then use the detector $\widehat{C}^{(0)}_{k,\ell}$ instead of $\widehat{C}_{k,\ell}$.
 But let us note that this situation is infrequent in practice.
 \end{enumerate}

  \end{rem}

 \noindent Therefore, at a nominal level $\alpha \in (0,1)$, take $c=c(\alpha)$ be the $(1-\alpha)$-quantile of the distribution of
  $ U_d = \sup_{0<u<1 } ~ f(u) \,  \| W_d(u)\|$ which can be computed through Monte-Carlo simulations.
  Table \ref{tab1}  shows the $(1-\alpha)$-quantile of this distribution for $ \alpha = 0.01, 0.05, 0.10$ and $d=1,\cdots,5$.

\begin{table}[htbp]
    \centering
    \begin{tabular}{c c c c c c}
     \hline
         & $d=1$  & $d=2$ &  $d=3 $ &  $d=4$ & $d=5$  \\
      \hline
      \hline
    $\alpha = 0.01$     &   $ 2.583$       &     $ 3.035 $      & $ 3.335 $ &  $ 3.631$ &  $ 3.914 $  \\
    &  &  &  & \\
    $\alpha = 0.05$    &   $ 1.954  $ &   $ 2.432 $  &   $ 2.760  $  &   $ 3.073  $  &  $ 3.334 $ \\
   &  &  &  & \\
    $\alpha = 0.10$    & $ 1.652  $ &  $ 2.156  $  &   $ 2.486 $  &   $ 2.784  $  &   $ 3.028 $  \\
  \hline
 \end{tabular}
 \caption{{\footnotesize  Empirical $(1-\alpha)$-quantile of the distribution of $ U_d $, for $d=1,\cdots,5$.  }}
        \label{tab1}
\end{table}

  \subsection{Asymptotic behaviour under the alternative hypothesis}
   Under the alternative $H_1$, the parameter  changes from $\theta^*_0$ to $\theta^*_1$ at $k^*>n$, where $\theta^*_1 \in \Theta$ and $\theta^*_0\neq \theta^*_1$. Then
  \begin{Theo}\label{theo2}
Assume \textbf{D}$(\Theta)$, \textbf{Id}($\Theta$),\textbf{Var}($\Theta$), \textbf{K}($f_\theta,M_\theta,\Theta$) and {\bf B}.
 Under the alternative $H_1$, if $\theta^*_1\neq \theta^*_0$ and $ \theta^*_0, \theta^*_1 \in \overset{\circ}{\Theta}(4)$ then
 for $ k^*=k^*(n)$  such as $\limsup_{n \rightarrow \infty} k^*(n)/n < \infty $ and
 $k_n=k^*(n)+n^\delta $ with $ \delta \in (1/2,1)$,
 $$  \max_{\ell \in \Pi_{n,k_n}} \, \dfrac{\widehat{C}_{k_n,\ell}}{b((k_n-\ell)/n)} \limitepsn  \infty  . $$
\end{Theo}
%
 The forthcoming Corollary \ref{cor2} can be immediately deduced from the relation (\ref{prob_stop_tim}).
 \begin{cor}\label{cor2}
  Under assumptions of Theorem \ref{theo2},
  $$  \lim_{n \rightarrow \infty} P\{\tau(n)<\infty \} = 1 .$$
 \end{cor}

\begin{rem}
 We know that the monitoring scheme rejects $H_0$ at the first time $k$ where
    $$ \max_{\ell \in \Pi_{n,k}} \dfrac{\widehat{C}_{k,\ell}}{b((k-\ell)/n)}>1 .$$
    Therefore, it follows from Theorem \ref{theo2} that under the hypothesis $H_1$, the detection delay $\widehat{d}_n$ of the procedure can be bounded by
    $\mathcal{O}_P(n^{1/2+\varepsilon})$ for any $\varepsilon>0$ (or even by $\mathcal{O}_P\big (\sqrt n (\log n)^a \big )$ with $a>0$ using the same kind of proof).
\end{rem}

\subsection{Examples}
   \subsubsection{AR($\infty$) processes}
          Consider the generalization of ARMA$(p,q)$ processes defined in \eqref{arma} {\it i.e.} a AR($\infty$) processes defined by:
          \begin{eqnarray}\label{arinfty}
           X_t = \phi_0(\theta^*_0) + \sum_{j\geq 1} \phi_j(\theta^*_0)X_{t-j} + \xi_t  ~ , ~ t \in \Z
           \end{eqnarray}
 with $ \theta^*_0 \in  \stackrel{\circ}{ \Theta}$, where we can chose $\Theta$ as a compact subset of $\Theta(4) \subset \R^d$ where
 $$\Theta(4)=\big \{ \theta\in \R^d; ~ \sum_{j\geq 1} |\phi_j(\theta)|  <1 \big \}.$$
  This process belongs to the class ${\cal M}_{\Z}(M_{\theta^*_0},f_{\theta^*_0})$ where
  $f_\theta(x_1,\cdots)= \sum_{j\geq 1} \phi_j(\theta)x_j $ and $ M_\theta \equiv \phi_0(\theta) ~ \text{for all}~ \theta \in \Theta$ and therefore  $\alpha^{(0)}_j(f_\theta,\Theta) = \|\phi_j(\theta)\|_\Theta$ and $\alpha^{(0)}_j(M_\theta,\Theta) =0$ for $j\in \N^*$. Then
  \begin{itemize}
 \item  Assumption {\bf D}($\Theta$) holds if $\underline{h}=\underset{\theta\in \Theta} {\mbox{inf}}(|\phi_0(\theta)|)>0$;
 \item  Assumption \textbf{K}($f_\theta,M_\theta,\Theta$) holds if there exists $\ell >2$ and and if $\theta \mapsto \phi_j(\theta)$ are twice differentiable functions satisfying  $\max\big ( \|\psi_j(\theta)\|_\Theta,\|\phi_j'(\theta)\|_\Theta, \|\phi_j''(\theta)\|_\Theta\big )  = O(j^{-\ell})$ for $j\in \N$.
 \item if $(\xi_t)$ is a sequence of non-degenerate random variables, Assumptions {\bf Id}($\Theta$) and {\bf Var}($\Theta$) hold.
   \end{itemize}

  ~\\
   \textbf{Case of AR($p$) process}\\
   Assume that
        $$ X_t= \phi^*_0+\sum\limits_{j=1}^{p}\phi^*_j X_{t-j} + \xi_t ~ ~ \text{with} ~ p \in \N^*   .$$
   The true parameter of
  the model is denoted by $ \theta^*_0 = (\phi^*_0,\phi^*_1, \cdots, \phi^*_{p}) \in \Theta$ where
  $ \Theta = \{  \theta = ( \phi_0,\phi_1, \cdots,  \phi_p) \in \R^{p+1} ~ ~ /  \sum\limits_{j=1}^{p}| \phi_j|<1   \}   $.
  Then, $ \Theta (r) = \Theta $ for any $r\geq 1$. Assume that a trajectory $(X_1, \cdots, X_k)$ has been observed,
   for any $t=1,\cdots,k$ and   $\theta \in \Theta$ we have,
  $ \widehat{q}_t(\theta)= \big(X_t- \phi_0-\sum\limits_{j=1}^{p}\phi_j X_{t-j} \big)^2$,
   $  \dfrac{\partial \widehat{q}_{t}( \theta)}{\partial \theta}  =
  -2 \big(X_t- \phi_0- \sum\limits_{j=1}^{p}\phi_j X_{t-j} \big)\cdot (1,X_{t-1}, X_{t-2},\cdots,X_{t-p})$.
   Moreover, $ \dfrac{\partial^2 \widehat{q}_{t}( \theta)}{\partial \phi_0 \partial \phi_0} = 2 $,
   for $  j=1,\cdots,p$, ~ $ \dfrac{\partial^2 \widehat{q}_{t}( \theta)}{\partial \phi_0 \partial \phi_j} = 2 X_{t-j} $
  and for $  1 \leq i,j\leq p$, ~  $ \dfrac{\partial^2 \widehat{q}_{t}( \theta)}{\partial \phi_i \partial \phi_j} = 2X_{t-i} X_{t-j}.$

  \subsubsection{ARCH($\infty$) processes}
          Consider the generalization  GARCH$(p,q)$ processes defined in \eqref{garch} {\it i.e.} a ARCH($\infty$) processes defined by:
          \begin{eqnarray}\label{archinfty}
           X_t = \sigma_t \, \xi_t \quad\mbox{and}\quad \sigma_t^2 = \psi_0(\theta^*_0) + \sum^{\infty}_{j= 1} \psi_j(\theta^*_0)X^2_{t-j}  ~ , ~ t \in \Z
           \end{eqnarray}
 with $ \theta^*_0 \in  \stackrel{\circ}{ \Theta}$, where we can chose $\Theta$ as a compact subset of $\Theta(4) \subset \R^d$ where
 $$\Theta(4)=\big \{ \theta\in \R^d; ~ (\E|\xi_0|^4)^{1/2} \sum_{j= 1}^\infty  |\phi_j(\theta)|  <1 \big \}.$$
  This process, introduced by Robinson \cite{Rob}, belongs to the class ${\cal M}_{\Z}(f_{\theta^*_0},M_{\theta^*_0},)$ where
  $f_\theta(x_1,\cdots)\equiv 0 $ and $ M^2_\theta(x_1,\cdots) = \psi_0(\theta) + \sum_{j\geq 1} \psi_j(\theta)x^2_j ~ \text{for all}~ \theta \in \Theta$ and therefore $\alpha^{(0)}_j(f_\theta,\Theta) = 0$ and $\alpha^{(0)}_j(h_\theta,\Theta) =\|\phi_j(\theta)\|_\Theta$ for $j\in \N^*$ ($X$ is of course a ARCH-type process). Then
  \begin{itemize}
 \item Assumption {\bf D}($\Theta$) holds if $\underline{h}=\underset{\theta\in \Theta} {\mbox{inf}}(\psi_0(\theta))>0$;
 \item Assumption \textbf{K}($f_\theta,M_\theta,\Theta$) holds if there exists $\ell >2$ and and if $\theta \mapsto \phi_j(\theta)$ are twice differentiable functions satisfying $ \max\big (\|\psi_j(\theta)\|_\Theta,\|\psi_j'(\theta)\|_\Theta ,\|\psi_j''(\theta)\|_\Theta\big )  = O(j^{-\ell})$ for $j\in \N$.
 \item if $(\xi^2_t)$ is a sequence of non-degenerate random variables, Assumptions {\bf Id}($\Theta$) and {\bf Var}($\Theta$) hold.
   \end{itemize}

   ~ \\
     \textbf{Case of GARCH($1,1$) process}\\
    Assume that
   $$  X_t= \sigma_t \xi_t ~~ \text{ with} ~~ \sigma^2_t= \alpha^*_0 + \alpha^*_1 X^2_{t-1} + \beta^*_1\sigma^2_{t-1}   $$
  with $ \theta^*_0=(\alpha^*_0,\alpha^*_1,\beta^*_1) \in \Theta \subset ]0,\infty[\times[0,\infty[^2$ and satisfying $ \alpha^*_1 + \beta^*_1 < 1 $.
  The ARCH($\infty$) representation is
    $  \sigma^2_t = \alpha^*_0/(1-\beta^*_1) + \alpha^*_1 \sum\limits_{ j\geq 1}^{} (\beta^*_1)^{j-1} X^2_{t-j}.  $
   If a trajectory $(X_1, \cdots, X_k)$ has been observed, for any $t=1,\cdots,k$ and   $\theta \in \Theta$ we have,
   $$ \widehat{h}^t_{\theta} = \alpha_0/(1-\beta_1) + \alpha_1X^2_{t-1}  + \alpha_1 \sum\limits_{ j= 2}^{t} \beta_1^{j-1} X^2_{t-j}  ~~ \text{and} ~~
      \widehat{q}_t(\theta) =  X^2_t /~ \widehat{h}^t_{\theta}  +  \log(\widehat{h}^t_{\theta})  .$$
    Thus, it follows that
   $ \dfrac{\partial \widehat{q}_t(\theta)}{ \partial \theta}  =  \dfrac{1}{\widehat{h}^t_{\theta} } \Big( 1 -\dfrac{X^2_t}{\widehat{h}^t_{\theta} } \Big)
      \Big( \dfrac{\partial \widehat{h}^t_{\theta}}{ \partial \alpha_0}  , \dfrac{\partial \widehat{h}^t_{\theta}}{ \partial \alpha_1}   ,
        \dfrac{\partial \widehat{h}^t_{\theta}}{ \partial \beta_1} \Big)   $
   with  $  \dfrac{\partial \widehat{h}^t_{\theta} }{  \partial \alpha_1 }= X^2_{t-1} + \sum\limits_{ j= 2}^{t} \beta_1^{j-1} X^2_{t-j} $
   $ \dfrac{ \partial \widehat{h}^t_{\theta}}{ \partial \alpha_0} = 1/(1-\beta_1) $,
   and ~ $\dfrac{ \partial \widehat{h}^t_{\theta}}{ \partial \beta_1} = \alpha_0/ (1-\beta_1)^2 + \alpha_1 X^2_{t-2} +  \alpha_1 \sum\limits_{ j= 3}^{t} (j-1)\beta_1^{j-2} X^2_{t-k} $.

    \noindent Let $ \theta = (\alpha_0,\alpha_1,\beta_1) = (\theta_1,\theta_2,\theta_3) \in \Theta  $, for $1 \leq i,j \leq 3$, we have
    $$ \dfrac{\partial^2 \widehat{q}_t(\theta)}{ \partial \theta_i  \partial \theta_j} =  \dfrac{1}{(\widehat{h}^t_{\theta})^2 }
  \Big( \dfrac{2X^2_t}{\widehat{h}^t_{\theta} } -1 \Big) \dfrac{\partial \widehat{h}^t_{\theta}}{ \partial \theta_i} \dfrac{\partial \widehat{h}^t_{\theta}}{ \partial \theta_j}
  +\dfrac{1}{\widehat{h}^t_{\theta} } \Big( 1 -\dfrac{X^2_t}{\widehat{h}^t_{\theta} } \Big) \dfrac{\partial^2 \widehat{h}^t_{\theta}}{\partial \theta_i \partial \theta_j } $$
   with ~ $\dfrac{ \partial^2 \widehat{h}^t_{\theta}}{ \partial \alpha_0^2}  = 0$,  $\dfrac{\partial^2 \widehat{h}^t_{\theta} }{ \partial \alpha_0  \partial \alpha_1 }=0$,    $\dfrac{\partial^2 \widehat{h}^t_{\theta} }{ \partial \alpha_1^2 } = 0 $,
     $ \dfrac{\partial^2 \widehat{h}^t_{\theta} }{ \partial \alpha_1  \partial \beta_1}  = X^2_{t-2} + \sum\limits_{ j= 3}^{t} (j-1)\beta_1^{j-2} X^2_{t-j} $,
    $ \dfrac{\partial^2 \widehat{h}^t_{\theta}}{\partial \alpha_0  \partial \beta_1}  = 1/(1-\beta_1)^2 $
     and ~ $\dfrac{ \partial \widehat{h}^t_{\theta} }{ \partial \beta^2_1} = 2 \alpha_0/ (1-\beta_1)^3 +2\alpha_1X^2_{t-3} +  \alpha_1 \sum\limits_{ j= 4}^{t} (j-1)(j-2)\beta_1^{j-3} X^2_{t-j} $.\\                          ~

   \subsubsection{TARCH($\infty$) processes}
 The process $X$ is called Threshold ARCH($\infty$) (TARCH($\infty$) in the sequel) if it satisfies
 \begin{equation} \label{tarch}
  X_t = \sigma_t \xi_t\quad\mbox{and}\quad \sigma_t = b_0(\theta_0^*) + \underset{j = 1}{\overset{\infty}{\sum}}
 \Big [b_j^+(\theta_0^*) \max ( X_{t - j},0)- b_j^-(\theta_0^*) \min (X_{t - j},0) \Big ]~,~t\in \Z
\end{equation}
 where the parameters $b_0(\theta)$, $b_j^+(\theta)$ and
 $b_j^-(\theta)$ are assumed to be non negative real numbers and $\theta \in  \stackrel{\circ}{ \Theta}$ where $\Theta$ is a compact subset of $\Theta(4)$ where
 $$
 \Theta(4)=\Big \{\theta\in\R^d~\Big/~ \big ( \E  | \xi_0|^4  \big )^{1/4} \, \sum_{j=1}^\infty \max \big (b_j^-(\theta),b_j^+(\theta)\big ) <1 \Big \}
 $$
 since $\alpha_j^{(0)}(M,\{\theta\})= \max \big (b_j^-(\theta),b_j^+(\theta)\big )$.
 This class of processes is a generalization of the class of
 TGARCH($p$,$q$) processes (introduced by Rabemananjara and Zako\"{\i}an
\cite{Rabemananjara1993}). Then,
 \begin{itemize}
 \item Assumption {\bf D}($\Theta$) holds if $\underline{h}=\inf_{\theta\in \Theta} b_0(\theta)>0$;
 \item Assumption \textbf{K}($f_\theta,M_\theta,\Theta$) holds if there exists $\ell >2$ and and if $\theta \mapsto b_j^-(\theta)$ and $\theta \mapsto b_j^+(\theta)$
 are twice differentiable functions satisfying
 $$ \max \big ( \|b_j^-(\theta)\|_\Theta,  \|b_j^+(\theta)\|_\Theta, \|\frac {\partial}{\partial \theta}  b_j^-(\theta)\|_\Theta,
  \|\frac {\partial}{\partial \theta}  b_j^+(\theta)\|_\Theta, \|\frac {\partial^2}{\partial \theta^2}  b_j^-(\theta)\|_\Theta,
  \|\frac {\partial^2}{\partial \theta^2}  b_j^+(\theta)\|_\Theta \big )  = O(j^{-\ell})\quad \mbox{for}\quad j\in \N.
 $$
   \end{itemize}
 Unfortunately, for TARCH($\infty$) it is not possible to provide simple conditions for obtaining   Assumptions {\bf Id}($\Theta$) and {\bf Var}($\Theta$) as for AR($\infty$) or ARCH($\infty$) processes.

\section{ Some simulation and numerical experiments}\label{simu}

  First remark that, at a time $k>n$, we need to compute $\widehat{C}_{k,\ell}$ for all $\ell \in \Pi_{n,k}$ to test whether change occurs or not.
  On can see that, the computational time is very long and increases with $k$.
  To reduce it, we introduce an integer sequence $(u_n)$ (satisfying $u_n/\sqrt{n} \rightarrow 0$ as $n\rightarrow \infty$; typically $u_n=[\ln (n)]$)
  and compute $\widehat{C}_{k,\ell}$ only for
  $$\ell \in \Pi^0_{n,k} := \{ n-v_n, n-v_n+u_n, n-v_n+2u_n, \cdots, k-v_n  \}.$$
  We have $\Pi^0_{n,k} \subset \Pi_{n,k}$ and for any $t=\frac{\ell}{n}$ with $\ell \in \Pi_{n,k}$, we can find an integer $j_{\ell}$  such that
  $n-v_n+ j_{\ell} u_n \in \Pi^0_{n,k} $  and  $ n-v_n+ j_{\ell} u_n \leq \ell \leq n-v_n+ (j_{\ell}+1) u_n $. This implies that
   $ \frac{n-v_n+ j_{\ell} u_n}{n} \leq t \leq \frac{n-v_n+ j_{\ell} u_n}{n} + \frac{u_n}{n}$. Thus, we have asymptotically (as $n\rightarrow \infty$),
   $t \sim \frac{n-v_n+ j_{\ell} u_n}{n}$.
  It shows that the previous asymptotic results still hold by computing $\widehat{C}_{k,\ell}$  for $\ell \in \Pi^0_{n,k}$.
  The condition $u_n/\sqrt{n} \rightarrow 0$ ensures that the Theorem \ref{theo2} still holds by choosing $k_n=k^*(n)+n^\delta $ with $ \delta \in (1/2,1)$.
  In practice, the use of $ \Pi^0_{n,k}$ can introduce a distortion in the detection delay.
  But, the new detection delay must be between $\widehat{d}_n$ and $\widehat{d}_n + u_n$ (where $\widehat{d}_n$ is the detection delay obtained by using $\Pi_{n,k}$).
  In the sequel, we use $u_n=[\ln (n)]$.\\
  Moreover, if $b\equiv c >0$ is a constant function, according to  (\ref{prob_stop_tim}), we have
   \begin{equation}
 P\{ \tau(n)<\infty \} = P\Big\{ ~ \underset{k>n} {\mbox{sup}} ~ \underset{\ell \in \Pi^0_{n,k}} {\mbox{max}} \widehat{C}_{k,\ell}>c \Big \}.
  \end{equation}
 Thus, denote
               $$ \widehat{C}_{k} = \underset{\ell \in \Pi^0_{n,k}} {\mbox{max}} \widehat{C}_{k,\ell} ~ ~ \text{for any} ~ ~ k>n.$$
   The procedure is monitored from $k=n+1$ to $k=n+500$. The set $\{n+1,\cdots,n+500 \}$ is called monitoring period.
   According to the Remark 1 of \cite{Kengne2011}, $v_n=[(\log n)^{\delta}]$ (with $1\leq \delta \leq 3$) is chosen.
  We evaluated the performance of the procedure with $ v_n=[\log n], ~ [(\log n)^{3/2}], ~ [(\log n)^2], ~ [(\log n)^3] $
  and we recommend to use $ v_n= [(\log n)^{3/2}]$
   for linear model and $ v_n= [(\log n)^2]$  for ARCH-type model. The nominal level used in the sequel is $\alpha =0.05$.

 \subsection{An illustration}

\begin{figure}
   \begin{center}
   \epsfig{file=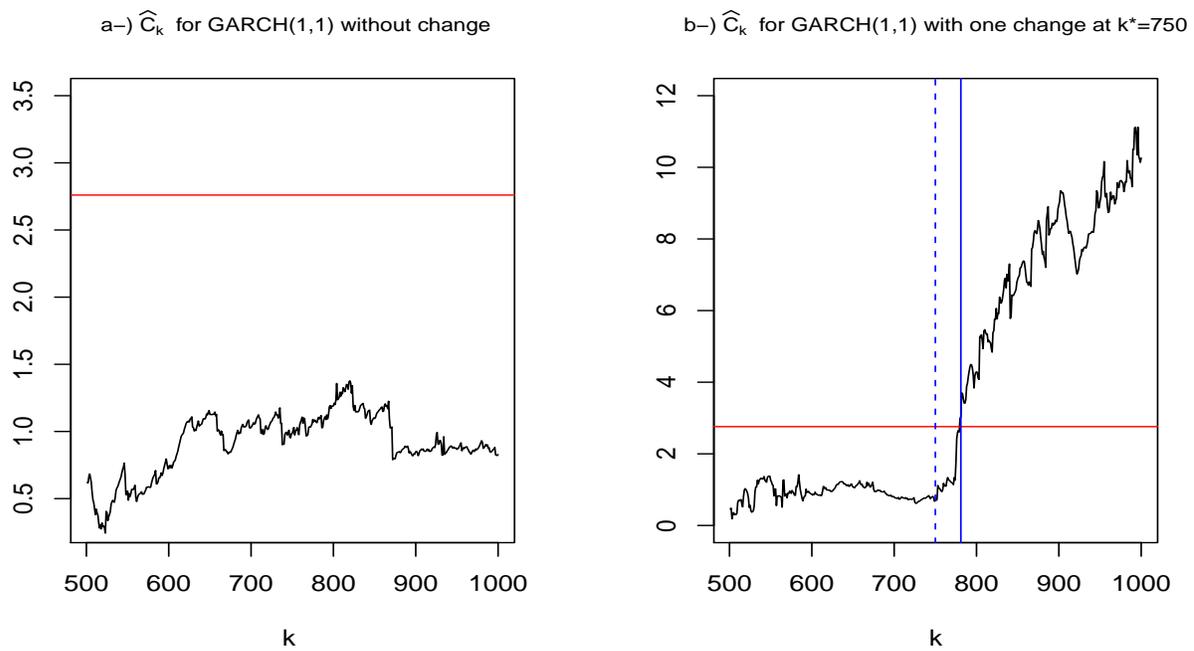, width=17.1cm, height=9.5cm}
   \end{center}
   \caption{Typical realization of the statistics  $\widehat{C}_{k}$  for GARCH(1,1), $n=500$ and $k=501,\cdots,1000$.
   a-) The parameter $\theta^*_0=(0.01,0.3,0.2)$ is constant ; b-) the parameter $\theta^*_0=(0.01,0.3,0.2)$  changes to $\theta^*_1 = (0.05,0.5,0.2)$ at $k^*=750$.
   The horizontal solid line represents the limit of the critical region, the vertical dotted line indicates where the change occurs and the   vertical solid
   line indicates the time where the monitoring procedure detecting the change. }
   \label{fig1}
   \end{figure}

  We consider a GARCH(1,1) process : $X_t= \sigma_t \xi_t ~~ \text{ with} ~~ \sigma^2_t= \alpha^*_0 + \alpha^*_1 X^2_{t-1} + \beta^*_1\sigma^2_{t-1}$.
   Thus, the parameter of the model is $\theta^*_0=(\alpha^*_0,\alpha^*_1,\beta^*_1)$.
  The historical available data are $X_1,\cdots,X_{500}$ (therefore $n=500$) and  the monitoring period is $\{501,\cdots,1000 \}$.
  At the nominal level $\alpha=0.05$, the critical values of the procedure is $C_{\alpha}=2.760$.
  The Figure \ref{fig1} is a typical realization of the statistic $(\widehat{C}_{k})_{500<k\leq 1000}$.
  We consider a scenario without change (Figure \ref{fig1} a-)) and a scenario with change at $k^*=n+250=750$ (Figure \ref{fig1}  b-)).\\
    Figure \ref{fig1} a-) shows that the detector $\widehat{C}_{k}$ is under the horizontal line which represents the limit of the critical region.
    On Figure \ref{fig1} b-) we can see that, before change occurs, $\widehat{C}_{k}$ is under the horizontal line and increases with a high speed after change.
    Such growth over a long period indicates that something happening in the model.

  \subsection{Monitoring mean shift in times series}
  Let $(X_1,\cdots,X_n)$ be an (historical) observation of a process $X=(X_t)_{t\in \Z}$. We assume that $X$ satisfy

   \begin{center}
     $\left\{
     \begin{array}{l}
      X_t = \mu_0 + \epsilon_t ~ ~ \text{for} ~ ~ 1\leq t \leq k^* \\
      X_t = \mu_1 + \epsilon_t ~ ~ \text{for} ~ ~ t > k^*
    \end{array}
    \right.$
 \end{center}
 with $k^*>n$, $\mu_0 \neq \mu_1$ and $(\epsilon_t)$ a zero mean stationary time series belongs to a class ${\cal M}_\Z (f_\theta,M_\theta)$.
 Under $H_0$,  $k^*=\infty$. The monitoring procedure start at $k=n+1$ and the aim is to test mean shift over the new observations $X_{n+1},X_{n+2},\cdots$.\\
~\\
 This problem can be seen as monitoring changes in linear model (see Horváth {\it et al.} \cite{Horvath2004}, Aue {\it et al.} \cite{Aue2006}) with constant regressor. The empirical mean $ \overline{X}_n =  \frac{1}{n}\sum_{i=1}^{n}X_i $
  is a consistent estimator of $\mu_0$ and the recursive residuals are defined by
  $ \widehat{\epsilon}_k = X_k- \overline{X}_n ~ ; ~ \text{for} ~ k>n.$
 Horváth {\it et al.} \cite{Horvath2004} and  Aue {\it et al.} \cite{Aue2006} proposed the CUSUM detector
\begin{equation}\label{cusum}
  \widehat{Q}_k = \dfrac{1}{\widehat{\sigma}_n } \dfrac{1}{c \sqrt{n} ( \frac{k}{n})(1-\frac{n}{k})^{\gamma} } \Big| \underset{i=n+1}{\overset{k} \sum} \widehat{\epsilon}_i \Big|
 \qquad k>n, ~ c>0, ~0\leq \gamma < 1/2,
\end{equation}
 where $\widehat{\sigma}^2_n$ is a consistent estimator of the long-run variance $ \sigma^2 = \lim_{n\rightarrow \infty } \frac{1}{n} \, \mbox{Var}( \underset{i=1}{\overset{n} \sum} \epsilon_i ).$
  If the process $(\epsilon_t)$ are uncorrelated (for instance GARCH-type model), empirical variance of the historical data can be used as estimator of $\sigma^2$.
  If $(\epsilon_t)$ are correlated, the popular Bartlett estimator (see \cite{Berkes2005}) can be used.
  Under some regular conditions, it hold that (see \cite{Horvath2004} and \cite{Aue2006})
  $$ \lim_{n\rightarrow \infty} P\{ \tau(n)<\infty \}=P\Big\{ ~ \underset{0<s<1} {\mbox{sup}} ~ \frac{| W_1(s) |}{s^\gamma}  >c \Big \}  . $$
   Hence, at a nominal level $\alpha=0.05$, the critical value of the test is the $(1-\alpha)$-quantile of the distribution of
  $\underset{0<s<1} {\mbox{sup}} ~ \| W_1(s) \|/s^{\gamma}$. When $\gamma=0$, these quantiles are known (see Table 1 of \cite{Lee2011}
  for values obtained through a Monte Carlo simulation).\\
~\\
  We compare our procedure to this CUSUM one (see \eqref{cusum} with $\gamma=0$) in two situations
  \begin{enumerate}
    \item  $(\epsilon_t)$ is an AR($1$) process;  $\epsilon_t= \phi^*_1 \epsilon_{t-1} + \xi_t$ with $\phi^*_1=0.2$;
    \item $(\epsilon_t)$ is a GARCH($1,1$) process;
    $\epsilon_t= \sigma_t \xi_t ~~ \text{ with} ~~ \sigma^2_t= \alpha^*_0 + \alpha^*_1 \epsilon^2_{t-1} + \beta^*_1\sigma^2_{t-1}$ and
    $(\alpha^*_0,\alpha^*_1,\beta^*_1)=(0.01,0.3,0.2)$.
  \end{enumerate}
  The historical sample size are $n=500$ and $n=1000$. These procedures are evaluated at times $k=n+100, n+200, n+300, n+400, n+500$, while the change occurs at $k^*=n+50$ or $k^*=n+250$.
  Tables \ref{tab2} and \ref{tab3} indicate the empirical levels and the empirical powers based of $200$ replications.
  The elementary statistics of the empirical detection delay are reported in Tables \ref{tab4}. \\
   \begin{table}[htbp]
    \centering
    \begin{tabular}{| c c c c c c c c|}
      \hline
           &   &      ~ ~ ~ $k$ &  $ n+100$ & $ n+200$  & $ n+300$ &  $ n+400$ & $ n+500$ \\
      \hline
       \hline
    Empirical levels  &  ~ $n=500$ ~ &  $\widehat{C}_k$  & 0.000   & 0.000   & 0.010    &   0.015   &  0.015   \\
                   &                &  $\widehat{Q}_k$   &  0.000  & 0.005  &  0.005     &  0.010   &  0.015     \\
                     &                  &                &    &   &       &      &    \\
                     &  ~ $n=1000$ ~ &  $\widehat{C}_k$   & 0.000   & 0.000   & 0.000   &  0.005    &  0.010   \\
                     &                  &  $\widehat{Q}_k$   & 0.000   &  0.000  &  0.005     &  0.000   &  0.010     \\
                     &                  &                &    &   &       &      &    \\
                  \hline
    Empirical powers  &  ~ $n=500$ ; $k^*=n+50$ ~ &  $\widehat{C}_k$  &  0.310  & 1   &   1  &   1   &  1   \\
                                 &                &  $\widehat{Q}_k$   & 0.335   & 1  &   1    &  1   &  1     \\
                               &                  &                &    &   &       &      &    \\
                      &  ~ $n=500$ ; $k^*=n+250$ ~ &  $\widehat{C}_k$   & 0.000   &  0.000  &  0.190  &   1   &   1  \\
                     &                  &  $\widehat{Q}_k$   &  0.000  &  0.000  &  0.130     & 0.965    &  1     \\
                   &                  &                &    &   &       &      &    \\
                    &  ~ $n=1000$ ; $k^*=n+50$ ~ &  $\widehat{C}_k$   & 0.075   & 1   & 1   &  1    &  1   \\
                     &                  &  $\widehat{Q}_k$   & 0.095   &  1  &   1    &  1   &   1    \\
                   &                  &                &    &   &       &      &    \\
                     &  ~ $n=1000$ ; $k^*=n+250$ ~ &  $\widehat{C}_k$   &  0.000  & 0.000   & 0.135   &  1    & 1    \\
                     &                  &  $\widehat{Q}_k$   &  0.000  & 0.000   &  0.075     &  0.980   &   1    \\
      \hline
    \end{tabular}
 \caption{{\footnotesize  Empirical levels and powers  for monitoring means shift in AR(1) with $\phi^*_1=0.2$. The empirical levels
                        are computed when $\mu_0=0$ and the empirical powers are computed when the mean  $\mu_0=0$ changes to  $\mu_1=1.2$. }}
        \label{tab2}
\end{table}

   \begin{table}[htbp]
    \centering
    \begin{tabular}{| c c c c c c c c|}
      \hline
           &   &      ~ ~ ~ $k$ &  $ n+100$ & $ n+200$  & $ n+300$ &  $ n+400$ & $ n+500$ \\
      \hline
       \hline
    Empirical levels  &  ~ $n=500$ ~ &  $\widehat{C}_k$  & 0.005  & 0.015    & 0.030   & 0.055  & 0.060    \\
                   &                &  $\widehat{Q}_k$   & 0.000  & 0.000   &  0.005      & 0.010     &  0.010      \\
                     &                  &                &    &   &       &      &    \\
                     &  ~ $n=1000$ ~ &  $\widehat{C}_k$   &  0.000   &  0.005   &  0.005   & 0.010    &  0.015   \\
                     &                  &  $\widehat{Q}_k$   &  0.000   & 0.000   & 0.000   &  0.010   & 0.010    \\
                     &                  &                &    &   &       &      &    \\
                  \hline
    Empirical powers  &  ~ $n=500$ ; $k^*=n+50$ ~ &  $\widehat{C}_k$  &  1   &  1   &  1    &  1   &  1   \\
                                 &                &  $\widehat{Q}_k$   & 1   &  1  &    1   &  1    &  1     \\
                               &                  &                &    &   &       &      &    \\
                      &  ~ $n=500$ ; $k^*=n+250$ ~ &  $\widehat{C}_k$   &  0.010   &  0.015   &  1   &    1   &   1  \\
                     &                  &  $\widehat{Q}_k$   &  0.000   &  0.000   & 0.920    &  1  & 1      \\
                   &                  &                &    &   &       &      &    \\
                    &  ~ $n=1000$ ; $k^*=n+50$ ~ &  $\widehat{C}_k$   & 0.995    &  1   &   1  &   1    &  1   \\
                     &                  &  $\widehat{Q}_k$   &  0.985   &   1  &     1   &  1   &  1      \\
                   &                  &                &    &   &       &      &    \\
                     &  ~ $n=1000$ ; $k^*=n+250$ ~ &  $\widehat{C}_k$   &  0.000   &  0.000   &  0.980  &  1     &   1   \\
                     &                  &  $\widehat{Q}_k$   &  0.000   &  0.000   &  0.765   & 1   &   1     \\
      \hline
    \end{tabular}
 \caption{{\footnotesize  Empirical levels and powers  for monitoring means shift in GARCH(1,1) with  $(\alpha^*_0,\alpha^*_1,\beta^*_1)=(0.01,0.3,0.2)$.
  The empirical levels are computed when $\mu_0=0$ and the empirical powers are computed when the mean  $\mu_0=0$ changes to  $\mu_1=0.3$. }}
        \label{tab3}
\end{table}

   \begin{table}[htbp]
    \centering
    \begin{tabular}{| c c c c c c c c c c|}
      \hline
          &                  &                &   &   &   &   &       &      &    \\
         $\widehat{d}_n$  &   &     & Mean  & SD  &   Min  &  $Q_1$  &  Med  &  $Q_3$ & Max \\
                     \hline
                      \hline
          AR(1)    & $n=500$ ;  $k^*=n+50$ ~ &  $\widehat{C}_k$  &  54.74  & 14.95  &  18  &  44   &  54    &  64   &  103   \\
                                 &                &  $\widehat{Q}_k$ & 53.78   &  14.72   &  16    &  43  &   54    &  63    &  102     \\
                               &                  &                &    &    &    &   &       &      &    \\
                     & $n=500$ ;  $k^*=n+250$ ~ &  $\widehat{C}_k$   & 63.14   & 23.18  &  12   & 45  &  61   &   77    &  135   \\
                     &                  &  $\widehat{Q}_k$   &   72.70  & 21.47  &  7   &  56   &  71.5  &  90   &  139      \\
                   &                  &                &   &   &   &   &       &      &    \\
                     & $n=1000$ ;  $k^*=n+50$ ~ &  $\widehat{C}_k$   &  75.84  &   14.19  &  37   &  66   &  75   &   83    &  114   \\
                     &                  &  $\widehat{Q}_k$   &  72.60  & 13.23  &  41   &  63   &   73   &  82    &  111      \\
                   &                  &                &   &   &   &   &       &      &    \\
                     & $n=1000$ ;  $k^*=n+250$ ~ &  $\widehat{C}_k$   & 76.24   &   19.15  &  23   &  60   &  76   &  89     &  140   \\
                     &                  &  $\widehat{Q}_k$   & 86.82  &  22.57  &  27    &  70   &  85    &  100    &   151     \\
                   &                  &                &   &   &   &   &       &      &    \\
                 \hline
          GARCH(1,1)    & $n=500$ ;  $k^*=n+50$ ~ &  $\widehat{C}_k$  &  20.21  & 6.15  &  1 & 16 & 20 & 24  & 35  \\
                                 &                &  $\widehat{Q}_k$ &  27.06  &  4.52  &  16  &  24  & 27 &  30  & 44    \\
                               &                  &                &    &    &    &   &       &      &    \\
                     & $n=500$ ;  $k^*=n+250$ ~ &  $\widehat{C}_k$   & 25.53   & 8.04  & 3  & 20  & 25  & 31   &  50   \\
                     &                  &  $\widehat{Q}_k$   & 35.40 &  10.01 &  13  &  28  &  35   &  41   & 62     \\
                   &                  &                &   &   &   &   &       &      &    \\
                     & $n=1000$ ;  $k^*=n+50$ ~ &  $\widehat{C}_k$   &  28.43   &  7.41  &  6   &  24   &  28   &  33     &  51   \\
                     &                  &  $\widehat{Q}_k$   &  36.98  &  5.09  &  21   &  33   &   37  &  40  &  48    \\
                   &                  &                &   &   &   &   &       &      &    \\
                     & $n=1000$ ;  $k^*=n+250$ ~ &  $\widehat{C}_k$   & 31.16  &  8.52  &  4   &  26   &  33   &  39   &  53   \\
                     &                  &  $\widehat{Q}_k$   &   44.35  &   10.04   &  14 &  37 & 45  &  50  &  71  \\
                   &                  &                &   &   &   &   &       &      &    \\
      \hline
    \end{tabular}
 \caption{{\footnotesize  Elementary statistics of the empirical detection delay for monitoring mean shift in AR(1) and GARCH(1,1).}}
        \label{tab4}
\end{table}

The results of Table \ref{tab2} and Table \ref{tab3} show that both the procedures based on detectors $\widehat{C}_k$ and $\widehat{Q}_k$
   are conservative. One can also see that the larger  $n$ (length of historical data) the smaller the distortion size  of these procedures.
   This is due to the fact that the length of monitoring period is fixed and does not increase with $n$.\\
   Under $H_1$, the change has been detected before the monitoring time $k=n+500$. But, as we mentioned above, the challenge of this problem is to minimize the detection delay. For this criteria, it can be seen in Table \ref{tab4} that in the case of the mean shift in AR process, our procedure works well as Horváth {\it et al.'s} procedure when the change occurs    at the beginning of the monitoring ($k^* = n+50$); but our procedure is a little more accurate when the change occurs a long time after the beginning of the  monitoring ($k^*=n+250$). For the case of the mean shift in GARCH process, our test procedure outperforms the Horváth {\it et al.'s} test in terms of mean and quantiles of the detection delay.

 \subsection{Monitoring parameter changes in AR(1) and GARCH(1,1) processes}

  In this subsection, we present some simulations results for monitoring parameter changes in  AR(1) and GARCH(1,1) models and compare our procedure
  to the one proposed by Na {\it et al.} \cite{Lee2011}.
  If the boundary function $b(\cdot) \equiv c >0$ with a real number $c>0$, Na {\it et al.} show that under H$_0$,
 $$ \lim_{n \rightarrow \infty} P\{ \tau(n)<\infty \}=  \lim_{n \rightarrow \infty} P\Big\{ ~ \underset{k>n} {\mbox{sup}} ~ \widehat{D}_k >c \Big \}
                     = P\Big\{ ~ \underset{0<s<1} {\mbox{sup}} ~ \| W_d(s) \| >c \Big \},$$
 where
  $$
    \widehat{D}_k:=\sqrt{n} \big\|\widehat{G}(T_{1,n})^{-1/2} \widehat{F}(T_{1,n})\big(\widehat{\theta}(T_{1,k}) - \widehat{\theta}(T_{1,n})\big)\big\| .
  $$
  Hence, at a nominal level $\alpha$, the critical value of their procedure is the $(1-\alpha)$-quantile of the distribution of
  $\underset{0<s<1} {\mbox{sup}} ~ \| W_d(s) \|$ which can be found in Table 1 of \cite{Lee2011}.\\
~\\
The comparisons between their procedure based on $\widehat{D}_k$ and ours based on $\widehat C_{k,\ell}$ are made in the following situations:
  \begin{enumerate}
 \item  For \textbf{AR($1$) model} :   $X_t= \phi^*_1 X_{t-1} + \xi_t$ . Under $H_0$,  $\theta_0=\phi^*_1=0.2$. Under $H_1$, $\theta_0$ changes to $\theta_1=-0.5$ at $k^*$.
    \item For \textbf{GARCH($1,1$) model} :
    $X_t= \sigma_t \xi_t ~~ \text{ with} ~~ \sigma^2_t= \alpha^*_0 + \alpha^*_1 X^2_{t-1} + \beta^*_1\sigma^2_{t-1}$.
     Under $H_0$,  $\theta_0=(\alpha^*_0,\alpha^*_1,\beta^*_1)=(0.01,0.3,0.2)$, while under $H_1$, $\theta_0=(0.01,0.3,0.2)$ changes to  $\theta_1=(0.05,0.5,0.2)$ at $k^*$.
  \end{enumerate}
  The sizes of historical samples  are $n=500$ and $n=1000$. The procedures are evaluated at times $k=n+100, n+200, n+300, n+400, n+500$, while the change occurs at $k^*=n+50$ or $k^*=n+250$.
  Tables \ref{tab5} and \ref{tab6} indicate the empirical levels and the empirical powers based of $200$ replications.
  The elementary statistics of the empirical detection delay are reported in Tables \ref{tab7}. \\
   \begin{table}[htbp]
    \centering
    \begin{tabular}{| c c c c c c c c|}
      \hline
           &   &      ~ ~ ~ $k$ &  $ n+100$ & $ n+200$  & $ n+300$ &  $ n+400$ & $ n+500$ \\
      \hline
       \hline
    Empirical levels  &  ~ $n=500$ ~ &  $\widehat{C}_k$  &  0.000   &  0.000   &   0.010   &   0.010    &  0.035    \\
                   &                &  $\widehat{D}_k$   & 0.000    &  0.000  &  0.000     &  0.000   &   0.025     \\
                     &                  &                &    &   &       &      &    \\
                     &  ~ $n=1000$ ~ &  $\widehat{C}_k$   & 0.000   &  0.000   &  0.000   &   0.010    &  0.025    \\
                     &                  &  $\widehat{D}_k$   & 0.000    &  0.000   &    0.000    &  0.000    &   0.020     \\
                     &                  &                &    &   &       &      &    \\
                  \hline
    Empirical powers  &  ~ $n=500$ ; $k^*=n+50$ ~ &  $\widehat{C}_k$  & 0.335    &  1   &  1   &  1   &  1   \\
                                 &                &  $\widehat{D}_k$   & 0.175    & 0.985   &   1     &   1   &   1     \\
                               &                  &                &    &   &       &      &    \\
                      &  ~ $n=500$ ; $k^*=n+250$ ~ &  $\widehat{C}_k$   &  0.000   &  0.000   &  0.180   &  990     &  1    \\
                     &                  &  $\widehat{D}_k$   &  0.000   &  0.000   &  0.095      &  0.865    &   1     \\
                   &                  &                &    &   &       &      &    \\
                    &  ~ $n=1000$ ; $k^*=n+50$ ~ &  $\widehat{C}_k$   &  0.065   &  0.995   &  1   &   1    &  1    \\
                     &                  &  $\widehat{D}_k$   &  0.090   &  0.975   &   1     &  1    &    1    \\
                   &                  &                &    &   &       &      &    \\
                     &  ~ $n=1000$ ; $k^*=n+250$ ~ &  $\widehat{C}_k$   &  0.000   &  0.000   &  0.140   &  0.990     &   1   \\
                     &                  &  $\widehat{D}_k$   &  0.000    &  0.000  &   0.075     &  0.855    &   0.995     \\
      \hline
    \end{tabular}
 \caption{{\footnotesize  Empirical levels and powers  for monitoring parameter change in AR(1) process. The empirical levels
                        are computed when $\theta_0=\phi^*_1=0.2$ is constant and the empirical powers are computed when   $\theta_0=0.2$ changes to $\theta_1=-0.5$. }}
        \label{tab5}
\end{table}

   \begin{table}[htbp]
    \centering
    \begin{tabular}{| c c c c c c c c|}
      \hline
           &   &      ~ ~ ~ $k$ &  $ n+100$ & $ n+200$  & $ n+300$ &  $ n+400$ & $ n+500$ \\
      \hline
       \hline
    Empirical levels  &  ~ $n=500$ ~ &  $\widehat{C}_k$  & 0.010    & 0.025    &   0.040   & 0.095    &  0.105   \\
                   &                &  $\widehat{D}_k$   & 0.010   &  0.015  &  0.040    & 0.040    &  0.055    \\
                     &                  &                &    &   &       &      &    \\
                     &  ~ $n=1000$ ~ &  $\widehat{C}_k$   & 0.000   & 0.000   & 0.030    & 0.045    &  0.055   \\
                     &                  &  $\widehat{D}_k$   &  0.000   & 0.000   &   0.010     &  0.015   &   0.035     \\
                     &                  &                &    &   &       &      &    \\
                  \hline
    Empirical powers  &  ~ $n=500$ ; $k^*=n+50$ ~ &  $\widehat{C}_k$  &  0.890   &   1  &   1   &  1   &  1   \\
                                 &                &  $\widehat{D}_k$   &  0.390  &  0.855  & 0.930    & 0.965    & 0.985    \\
                               &                  &                &    &   &       &      &    \\
                      &  ~ $n=500$ ; $k^*=n+250$ ~ &  $\widehat{C}_k$   &  0.010   &  0.030   &  0.825   &  1    &  1   \\
                     &                  &  $\widehat{D}_k$   & 0.010  & 0.020   & 0.270   & 0.805  &  0.915    \\
                   &                  &                &    &   &       &      &    \\
                    &  ~ $n=1000$ ; $k^*=n+50$ ~ &  $\widehat{C}_k$   &  0.835   &  1   &  1   &  1   &  1   \\
                     &                  &  $\widehat{D}_k$   &  0.310  &  0.970  &  0.990   &  0.995   & 0.995     \\
                   &                  &                &    &   &       &      &    \\
                     &  ~ $n=1000$ ; $k^*=n+250$ ~ &  $\widehat{C}_k$   &  0.000   & 0.005    &  0.685   &  1     &  1    \\
                     &                  &  $\widehat{D}_k$   &  0.000   &  0.000   & 0.250   &   0.955   &  0.990    \\
      \hline
    \end{tabular}
 \caption{{\footnotesize   Empirical levels and powers  for monitoring parameter change in GARCH(1,1) process. The empirical levels
                        are computed when  $\theta_0=(\alpha^*_0,\alpha^*_1,\beta^*_1)=(0.01,0.3,0.2)$ is constant (hypothesis $H_0$) and the empirical powers are computed when
                          $\theta_0=(0.01,0.3,0.2)$  changes to $\theta_1=(0.05,0.5,0.2)$ (hypothesis $H_1$).}}
        \label{tab6}
\end{table}

   \begin{table}[htbp]
    \centering
    \begin{tabular}{| c c c c c c c c c c|}
      \hline
    &                  &                &   &   &   &   &       &      &    \\
          $\widehat{d}_n$   &   &     & Mean  & SD  &   Min  &  $Q_1$  &  Med  &  $Q_3$ & Max \\
                     \hline
                      \hline
          AR(1)    & $n=500$ ;  $k^*=n+50$ ~ &  $\widehat{C}_k$  &  55.36   & 18.75   &  9   &  42    &   56   &  67    &  121    \\
                                 &                &  $\widehat{D}_k$ &  71.54   &  38.44    &   2   &  52.75   &  69     &  89    &   167     \\
                               &                  &                &    &    &    &   &       &      &    \\
                     & $n=500$ ;  $k^*=n+250$ ~ &  $\widehat{C}_k$   &  66.81   & 25.27   &   5   &  49  &  65    &   83     &  149    \\
                     &                  &  $\widehat{D}_k$   &   97.80    &  39.42  &  21   &  68    & 89   &  123    &   222      \\
                   &                  &                &   &   &   &   &       &      &    \\
                     & $n=1000$ ;  $k^*=n+50$ ~ &  $\widehat{C}_k$   &  75.13   &  19.87    &  24    &  62   &  74    &   90     &  147    \\
                     &                  &  $\widehat{D}_k$   &  87.70   &  28.72   &   14    &  66    &  85     &   109    &   195      \\
                   &                  &                &   &   &   &   &       &      &    \\
                     & $n=1000$ ;  $k^*=n+250$ ~ &  $\widehat{C}_k$   &  76.89   &  26.16    &  15    &   56   &   77   &   96     &  172    \\
                     &                  &  $\widehat{D}_k$   &  101.20  &  37.97   &  20   &  75    &   96    &  129     &   245      \\
                   &                  &                &   &   &   &   &       &      &    \\
                 \hline
          GARCH(1,1)    & $n=500$ ;  $k^*=n+50$ ~ &  $\widehat{C}_k$  &  29.41  &  15.84  &  4  &  22   &  31    &  40   &  98   \\
                                 &                &  $\widehat{D}_k$ & 86.05  & 90.50   &  2  &  36  & 61    &  99    &  416    \\
                               &                  &                &    &    &    &   &       &      &    \\
                     & $n=500$ ;  $k^*=n+250$ ~ &  $\widehat{C}_k$   & 38.02   & 19.33  &  5   &  27   &  37   &   44    &  113   \\
                     &                  &  $\widehat{D}_k$   & 87.72  &  50.96  &  1   &  49.25   &  79   &  112    &  236    \\
                   &                  &                &   &   &   &   &       &      &    \\
                     & $n=1000$ ;  $k^*=n+50$ ~ &  $\widehat{C}_k$   &  41.96  &  13.93  &  3   &  32   &  41   &  48     &  94   \\
                     &                  &  $\widehat{D}_k$   & 71.29  &  37.12  &  6   &  46    &    66    &  88    &  287      \\
                   &                  &                &   &   &   &   &       &      &    \\
                     & $n=1000$ ;  $k^*=n+250$ ~ &  $\widehat{C}_k$   &  44.99  & 17.16  &  5   &  35   &  41   &  52   &  117   \\
                     &                  &  $\widehat{D}_k$   & 75.78  &  35.10   &  7  &  52  & 71   &   95   &    198    \\
                   &                  &                &   &   &   &   &       &      &    \\
      \hline
    \end{tabular}
 \caption{{\footnotesize  Elementary statistics of the empirical detection delay for monitoring parameter change in AR(1) and GARCH(1,1).}}
        \label{tab7}
\end{table}
  ~ \\
   \indent The considered  AR and GARCH processes  have zero mean. Contrary to the mean shift studied above, this mean is not estimated.
   For AR model, it appears in Table \ref{tab5} that both procedures based on detector $\widehat{C}_k$ and $\widehat{D}_k$ are conservative.
   This is not the case for GARCH model (Table \ref{tab6}). The high size distortions when $n=500$ is due to the difficulty to estimate the parameter of GARCH model.
   This size distortion decreases when $n$ increases and Corollary \ref{cor1} ensures that with infinite monitoring period, the empirical
   level tends to the nominal one as $n \rightarrow \infty$.\\
   \indent For both the cases of AR and GARCH processes, the procedure based on detector $\widehat{C}_{k,\ell}$ detects the change before the monitoring time $k=n+500$.
   Unlike Na {\it et al.} \cite{Lee2011}, we consider a scenario of GARCH model with moderate change in parameter and it can be seen in Table \ref{tab6} that the procedure based   on detector $\widehat{D}_k$ provides unsatisfactory results. At the monitoring time $k=n+500$, it is not sure that the change must be detected even when $k^*=n+50$.
   This is not surprising according to the comment of subsection \ref{monit_proc}.\\
   \indent Table \ref{tab7} indicates the distribution of the detection delay $\widehat{d}_n$. We can see in  Table \ref{tab7} (even in Table \ref{tab4}) that
   for our procedure, the relation $\widehat{d}_{1000}\leq \sqrt{1000/500}~ \widehat{d}_{500} $ is globally satisfied (from Theorem \ref{theo2}, we deduced that $\widehat{d}_{n}= {\cal O}_P\big (n^{1/2}\log n\big )$ when $n$ is large enough). Moreover, elementary statistics (mean and quantiles) show that the detection delay using our procedure  is shorter than using  the Na {\it et al.}'s one. The results of Table \ref{tab5}, \ref{tab6} and \ref{tab7} show that, our test is uniformly better and the procedure based on detector $\widehat{C}_k$ could be recommended in this frame.

    \section{Real-Data Applications }
\begin{figure}
   \begin{center}
   \epsfig{file=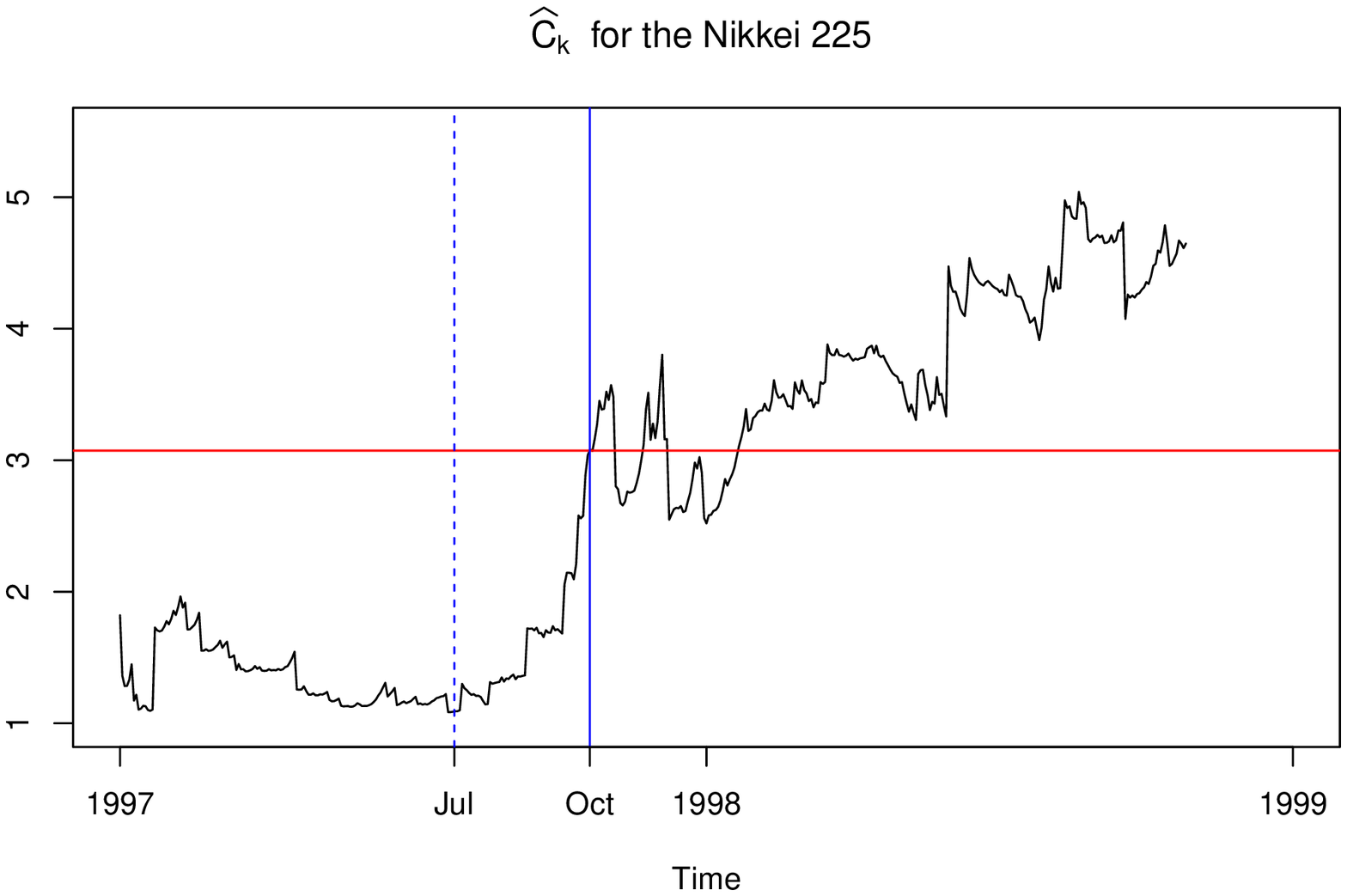, width=13.1cm, height=8.5cm}\\
   \epsfig{file=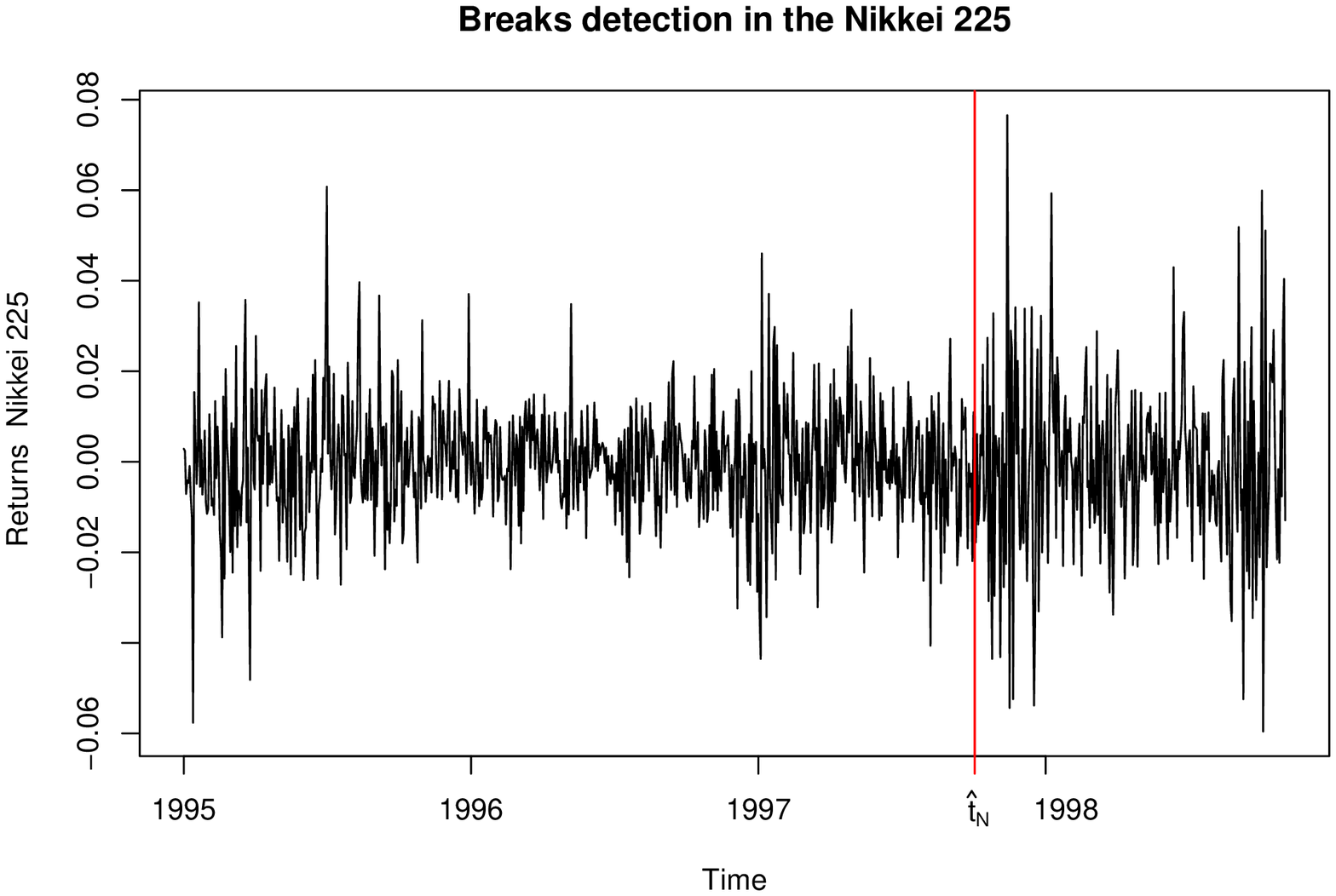 ,  width=13.1cm, height=8.5cm}
   \end{center}
   \caption{The top figure is a realization of the statistics  $\widehat{C}_{k}$  with $k$  going from January 2, 1995 to October 19, 1998 for Nikkei 225 data;
   the historical data considered are the series going from January 2, 1995 to December 31, 1996. The horizontal solid line represents the limit of the critical
   region, the vertical dotted line indicates the date of the beginning of the Asian financial crisis (1997-1998) and the  vertical solid line indicates the
   time where the monitoring procedure will stop. The bottom figure is the returns of Nikkei 225 data from January 2, 1995 to October 19, 1998; the vertical solid line
   indicates the date where break have been detected using retrospective test after the monitoring stops.}
   \label{fig2}
   \end{figure}

   We consider the returns of the daily closing values of the Nikkei 225  stock index (from January 2, 1995 to October 19, 1998),
   S$\&$P 500 and FTSE 100 (from January 2, 2004 to June 11, 2012). These data are available on Yahoo! Finance at {\tt http://finance.yahoo.com/}.
   They are represented on Figure \ref{fig2} and Figure \ref{fig6}. These series are known to represent ARCH effect and GARCH(1,2) (resp. GARCH(1,1)) can be used to capture it
   in returns of Nikkei 225 (resp.  S$\&$P 500 and FTSE 100), see the book of Francq and Zako\"{i}an 2010. \\
Consider the observations going from January 2, 1995 to December 31, 1996 (resp. January 2, 2004 to December 30, 2005)
   as the historical data for the Nikkei 225 (resp. S$\&$P 500 and FTSE 100) stock index. These periods are known to be stable in the financial community.
   To verify it, we apply three procedures to test for parameter change in the historical observations.
   The null hypothesis is that the parameter is constant over the observations against the parameter changes alternative.
    \begin{itemize}
        \item The first test is proposed by Kengne \cite{Kengne2011}. Define the asymptotic covariance matrix (which take into account  the change possibility)
         of the estimator $\widehat{\theta}_n(T_{1,n})$ by
              $$  \widehat{\Sigma}_{n,k} := \dfrac{k}{n}\widehat{F}_n(T_{1,k})\widehat{G}_n(T_{1,k})^{-1}\widehat{F}_n(T_{1,k})\text{\large{\1}}_{\det(\widehat{G}_n(T_{1,k}))\neq 0 }
                 +  \dfrac{n-k}{n}\widehat{F}_n(T_{k,n})\widehat{G}_n(T_{k,n})^{-1}\widehat{F}_n(T_{k,n})\text{\large{\1}}_{\det(\widehat{G}_n(T_{k,n}))\neq 0 } . $$
              The test is based on the statistic
              $$ \widehat{Q}_n:=\text{max} \big(\widehat{Q}^{(1)}_n , \widehat{Q}^{(2)}_n \big) ~ ~  \text{where}$$
 $$  \widehat{Q}^{(1)}_n:= \underset{ v_n \leq k \leq n-v_n } {\mbox{max}}\widehat{Q}^{(1)}_{n,k} ~~ ~
 \text{with} ~~  \widehat{Q}^{(1)}_{n,k}:= \dfrac{k^2}{n} \big( \widehat{\theta}_n(T_{1,k})-  \widehat{\theta}_n(T_n) \big)' \widehat{\Sigma}_{n,k} \big( \widehat{\theta}_n(T_k)-  \widehat{\theta}_n(T_n) \big),$$
   $$  \widehat{Q}^{(2)}_n:= \underset{ v_n \leq k \leq n-v_n } {\mbox{max}}\widehat{Q}^{(2)}_{n,k} ~~ ~
 \text{with} ~~  \widehat{Q}^{(2)}_{n,k}:= \dfrac{(n-k)^2}{n} \big( \widehat{\theta}_n(T_{k,n})-  \widehat{\theta}_n(T_{1,n}) \big)' \widehat{\Sigma}_{n,k} \big( \widehat{\theta}_n(T_{k,n})-  \widehat{\theta}_n(T_{k,n}) \big).$$
    This test is applied with $v_n=(\log n)^{\delta}$ where $2\leq \delta \leq 5/2$.
        \item The second test (see Lee and Song \cite{Lee2008}) is based on the statistic
          $$ \widehat{Q}^{(0)}_n:= \underset{ v_n \leq k \leq n-v_n } {\mbox{max}} \Big(   \dfrac{k^2}{n} \big( \widehat{\theta}_n(T_{1,k})-  \widehat{\theta}_n(T_{1,n}) \big)'
  \widehat{\Sigma}_{n,n} \big( \widehat{\theta}_n(T_{1,k})-  \widehat{\theta}_n(T_{1,n}) \big) \Big) $$
       with $v_n=(\log n)^2$.
        \item The third procedure is the CUSUM test see Kulperger and Yu \cite{Kulperger2005}.
    \end{itemize}
   At a nominal level $\alpha \in (0,1)$, each of these procedure rejects null hypothesis if the test statistic is greater than a critical
   value $C_{\alpha}$. Table \ref{tab8} provides the results of these tests to the historical data that we have chosen.

   \begin{table}[htbp]
    \centering
    \begin{tabular}{| c c c c |}
      \hline
         &  $ \widehat{Q}_n $  &  $\widehat{Q}^{0}_n$   &  CUSUM     \\
                     \hline
                      \hline
         Nikkei 225   & 3.35 (3.98)  & 2.31 (3.45) & 0.98 (1.36)  \\
         S$\&$P 500   &  2.13 (3.47)   & 2.01 (3.06)  &   0.93 (1.36)  \\
         FTSE 100  & 1.95 (3.47)  &   2.31 (3.06)   &  1.13 (1.36)    \\
      \hline
    \end{tabular}
 \caption{{\footnotesize Results of test for parameter changes in the historical data of Nikkei 225 (from January 2, 1995 to December 31, 1996),
             S$\&$P 500 and FTSE 100 (from January 2, 2004 to December 30, 2005). Figures in brackets the critical values of the procedure at the
             nominal level $\alpha = 0.05.$}}
        \label{tab8}
\end{table}
 ~ \\
 Note that, these series are very closed to a nonstationary process, in the sense that $ \sum_{j=1}^{q}\alpha_j +  \sum_{j=1}^{p}\beta_j\simeq 1$ (see (\ref{garch})). Therefore, it would be difficult to compute the estimator $\widehat{\theta}_n(T_{k,l})$ ($1\leq k < l \leq n$).
 For the statistics $\widehat{Q}_n$ and $\widehat{Q}^{(0)}_n$, we consider only the time $k$ that the computation of $\widehat{\theta}_n(T_{1,k})$
 and $\widehat{\theta}_n(T_{k,n})$ converges. This certainly introduces distortions on these tests. On the other hand, the CUSUM procedure needs to compute only
 the estimator $\widehat{\theta}_n(T_{1,n})$ which convergence is obtained. According to  these results, we conclude that the parameter
  does not change over these historical observations.\\
  \indent For Nikkei 225 data, monitoring starting at January 2, 1997. Figure \ref{fig2} shows the realization of the sequence $(\widehat{C}_k)$ for $k$
  going from January 2, 1997 to October 19, 1998. Monitoring procedure stops at October 27, 1997.
    Recall that, the monitoring scheme can be used as an alarm system. When it triggered, we need to apply  retrospective test to estimate the breakpoint.
   According to Kengne \cite{Kengne2011}, the test based on $\widehat{Q}_n$ and $\widehat{Q}^{(0)}_n$ are more powerful than the CUSUM test.
   Thus, in the retrospective procedure, we applied these two tests and considered the one which provides more significant result (in terms of p-value).\\
   Retrospective procedure is applied to the observations going from   January 2, 1995 to October 27, 1997 and break is detected at
   $\widehat{t}_{N} \simeq $ September 17, 1997;
   see Figure \ref{fig2}. This change corresponds to the Asian financial crisis (1997-1998) where the turmoil period started at July 1997.

\begin{figure}
   \begin{center}
   \epsfig{file=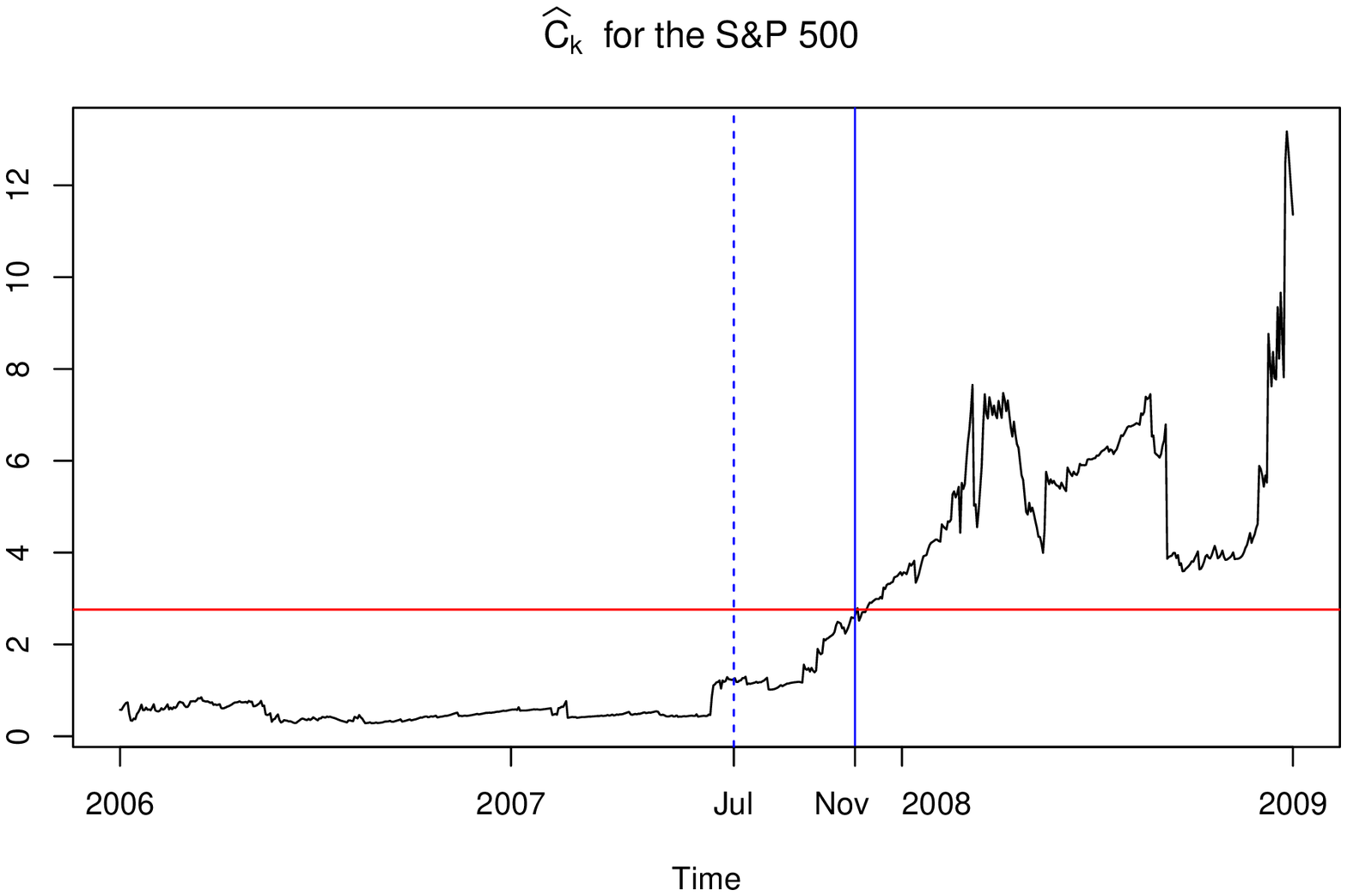, width=13.1cm, height=7.5cm}\\
   \epsfig{file=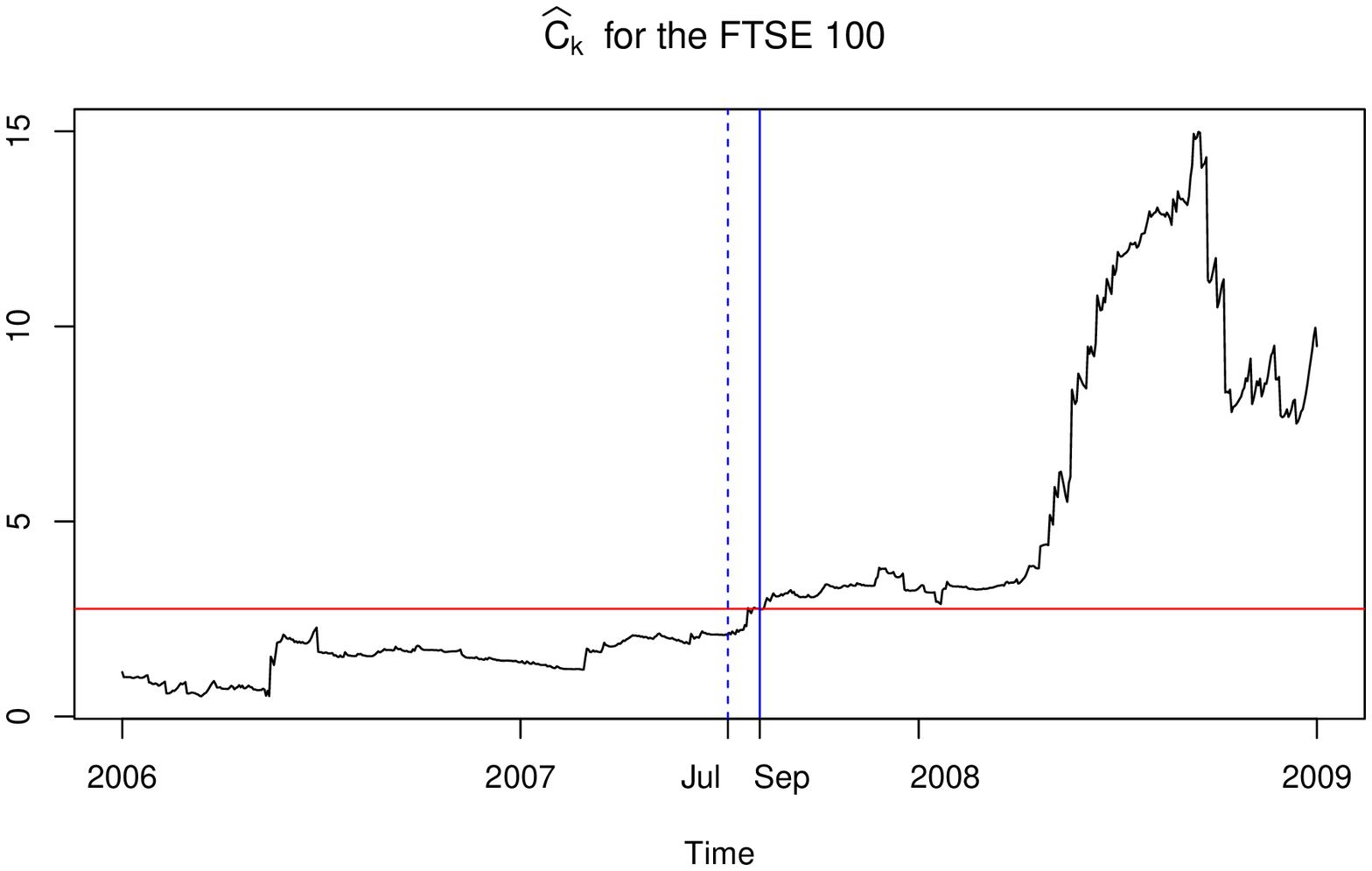 ,  width=13.1cm, height=7.5cm}
   \end{center}
   \caption{  Realization of the statistics  $\widehat{C}_{k}$  with $k$  going from January 2, 2006 to December 31, 2008 for S$\&$P 500 and FTSE 100 data;
    the historical data considered are the series going from January 2, 2004 to December 30, 2005.  The horizontal solid line represents the limit of
    the critical region, the vertical dotted line indicates the date of the beginning of the Subprime Crisis in US and the  vertical solid line indicates the
    time where the monitoring procedure stopped.}
   \label{fig3}
   \end{figure}

\begin{figure}
   \begin{center}
   \epsfig{file=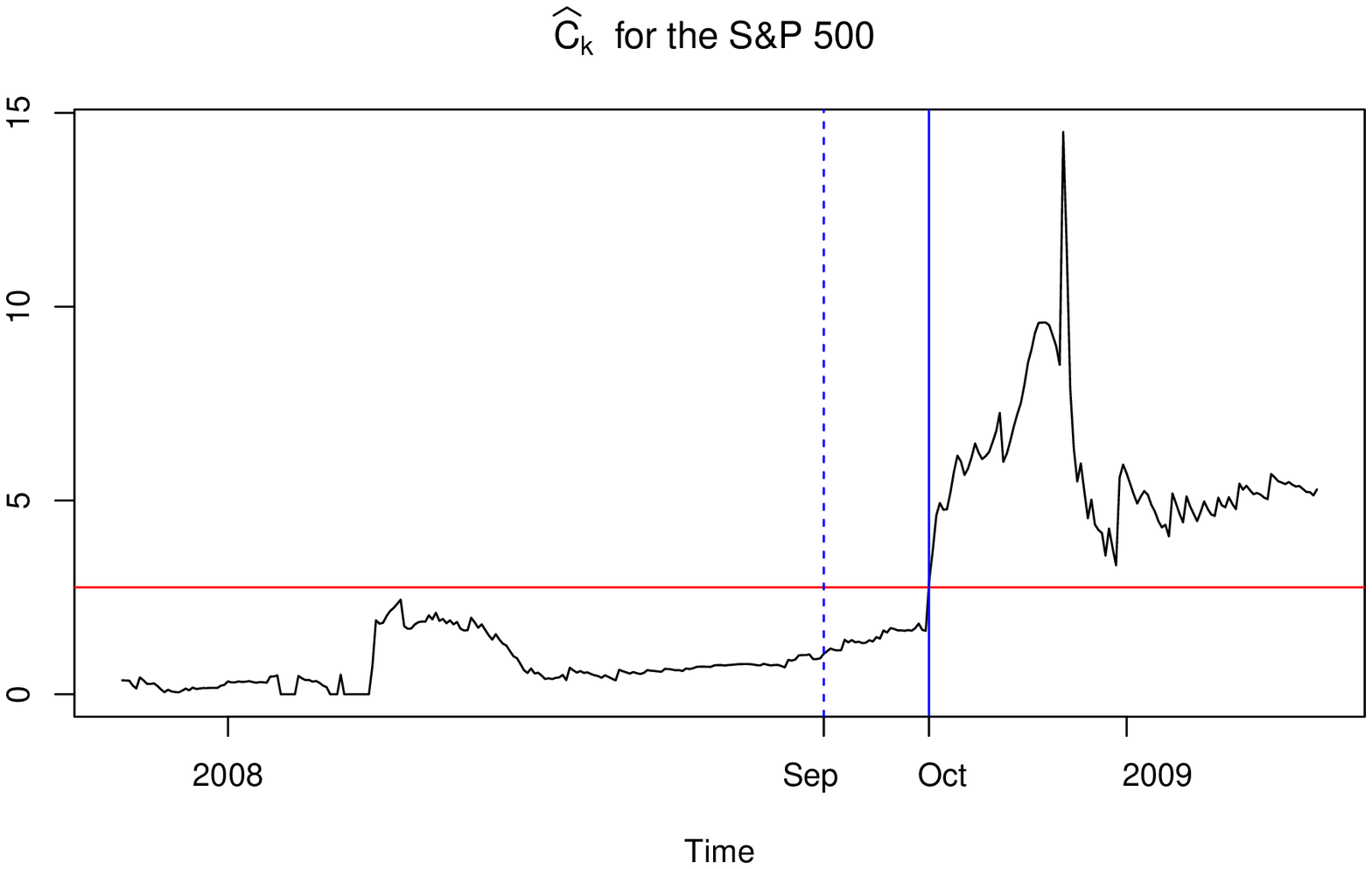, width=13.1cm, height=7.5cm}\\
   \epsfig{file=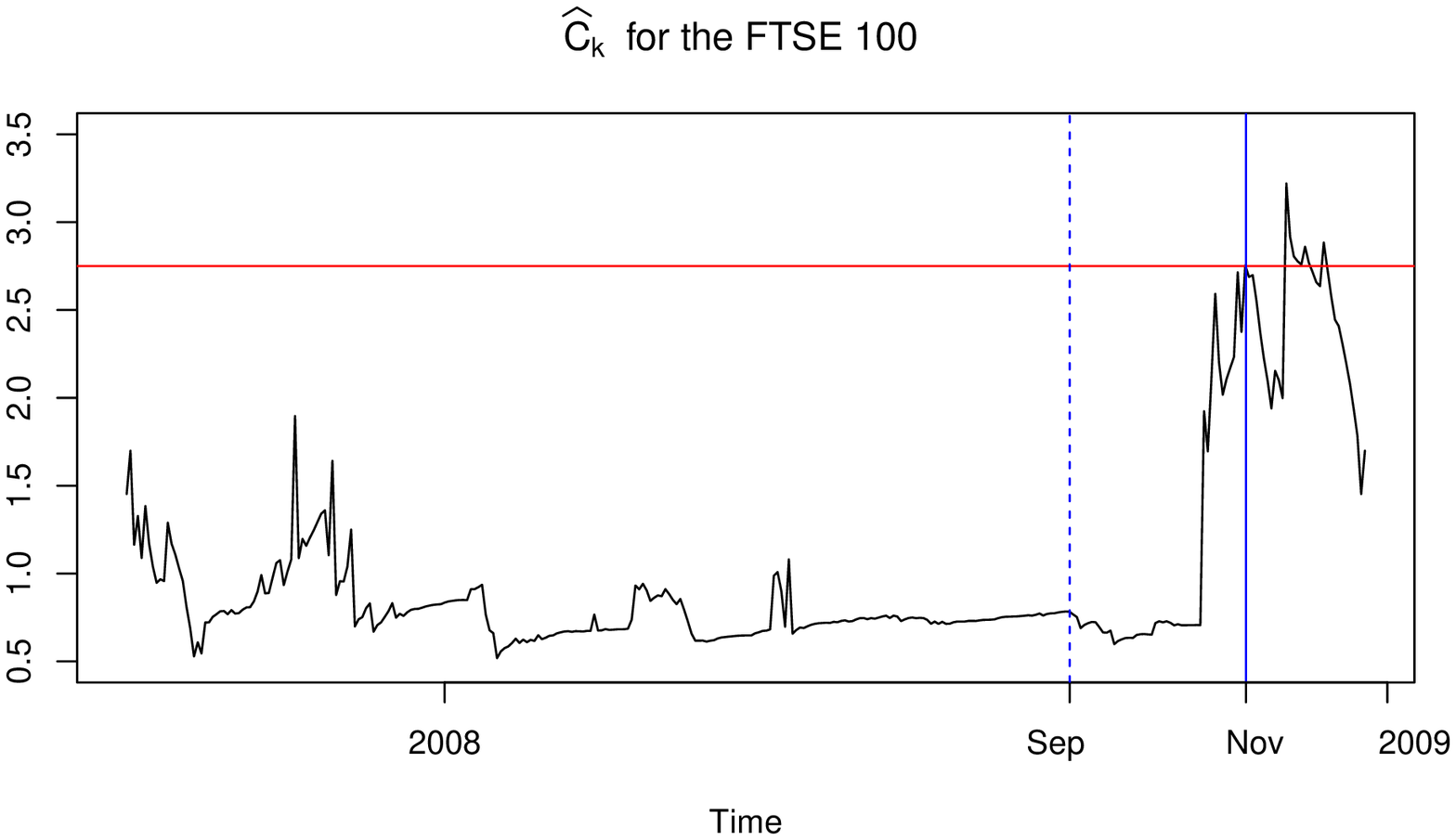 ,  width=13.1cm, height=7.5cm}
   \end{center}
   \caption{  Realization of the statistics  $\widehat{C}_{k}$ for S$\&$P 500 and FTSE 100 data;
   the historical data are the series going from June 18, 2007 to November 16, 2007 (for S$\&$P 500 data) and July 6, 2007 to September 4, 2007 (for FTSE 100 data).
   The horizontal solid line represents the limit of the critical region, the vertical dotted line indicates the date of the
   Lehman Brothers Crisis and the  vertical solid line indicates the time where the monitoring procedure stopped.  }
   \label{fig4}
   \end{figure}

\begin{figure}
   \begin{center}
   \epsfig{file=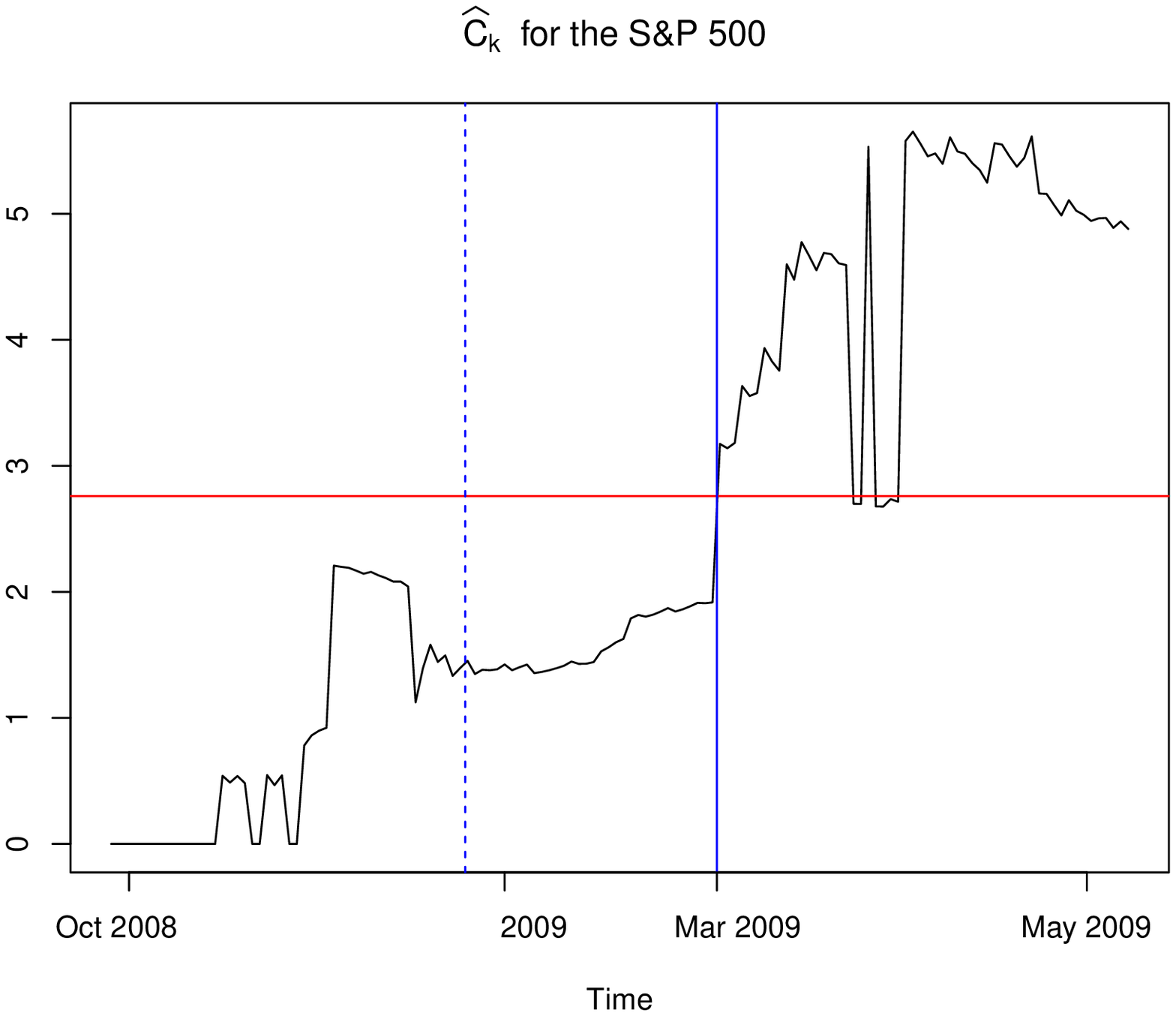 ,  width=8.1cm, height=8.5cm}
   \epsfig{file=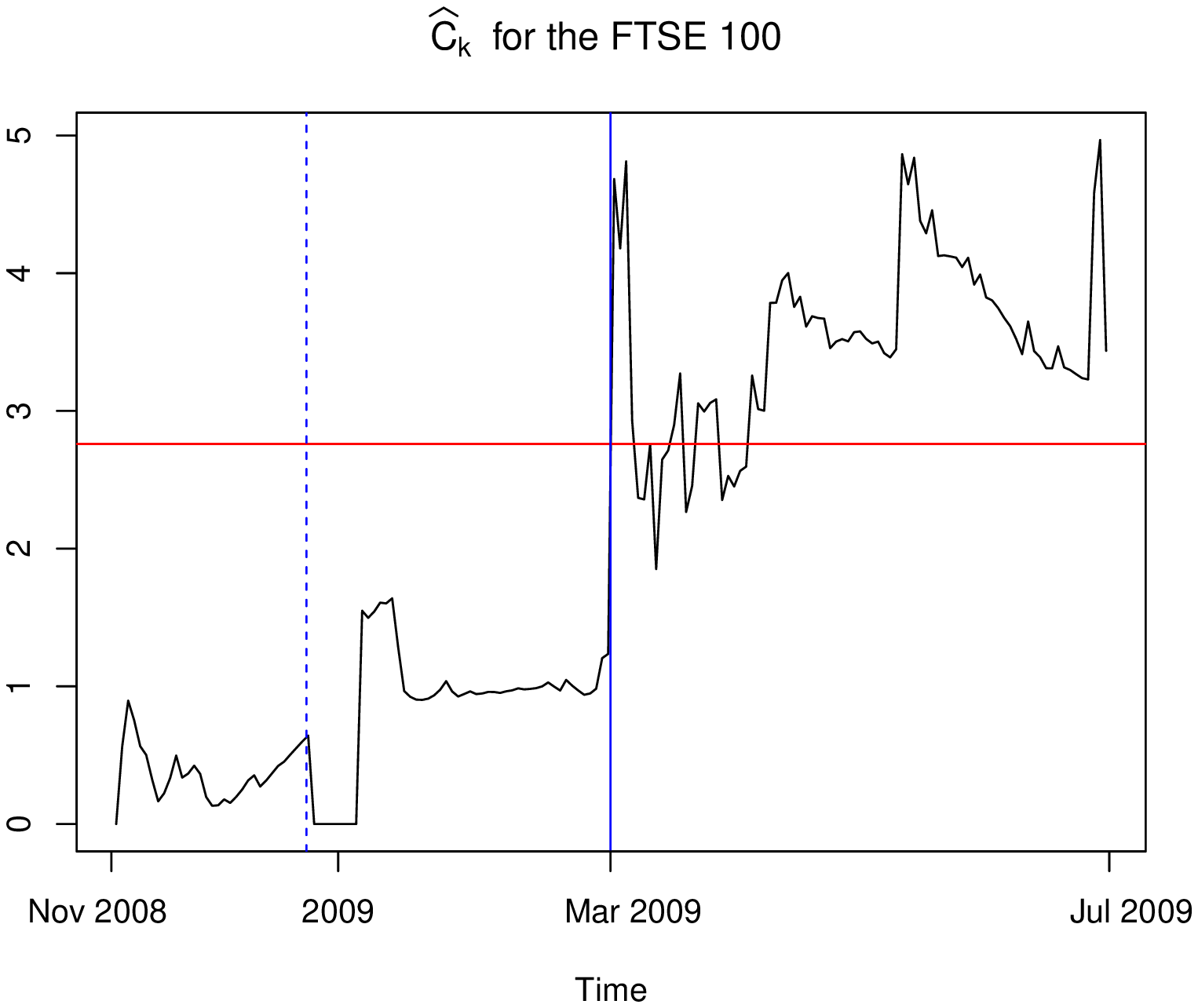, width=8.1cm, height=8.5cm}
   \end{center}
   \caption{ Realization of the statistics  $\widehat{C}_{k}$ for S$\&$P 500 and FTSE 100 data; the historical data are the series going
   from August 14, 2008 to October 17, 2008 (S$\&$P 500) and September 17, 2008 to November 10, 2008 (FTSE 100).
    The horizontal solid line represents the limit of the critical region, the vertical dotted line indicates the date of
    the beginning of stabilization in financial system due to the  worldwide  governments intervention and the vertical solid line indicates
    the time where the monitoring procedure stopped.}
   \label{fig5}
   \end{figure}

\begin{figure}
   \begin{center}
   \epsfig{file=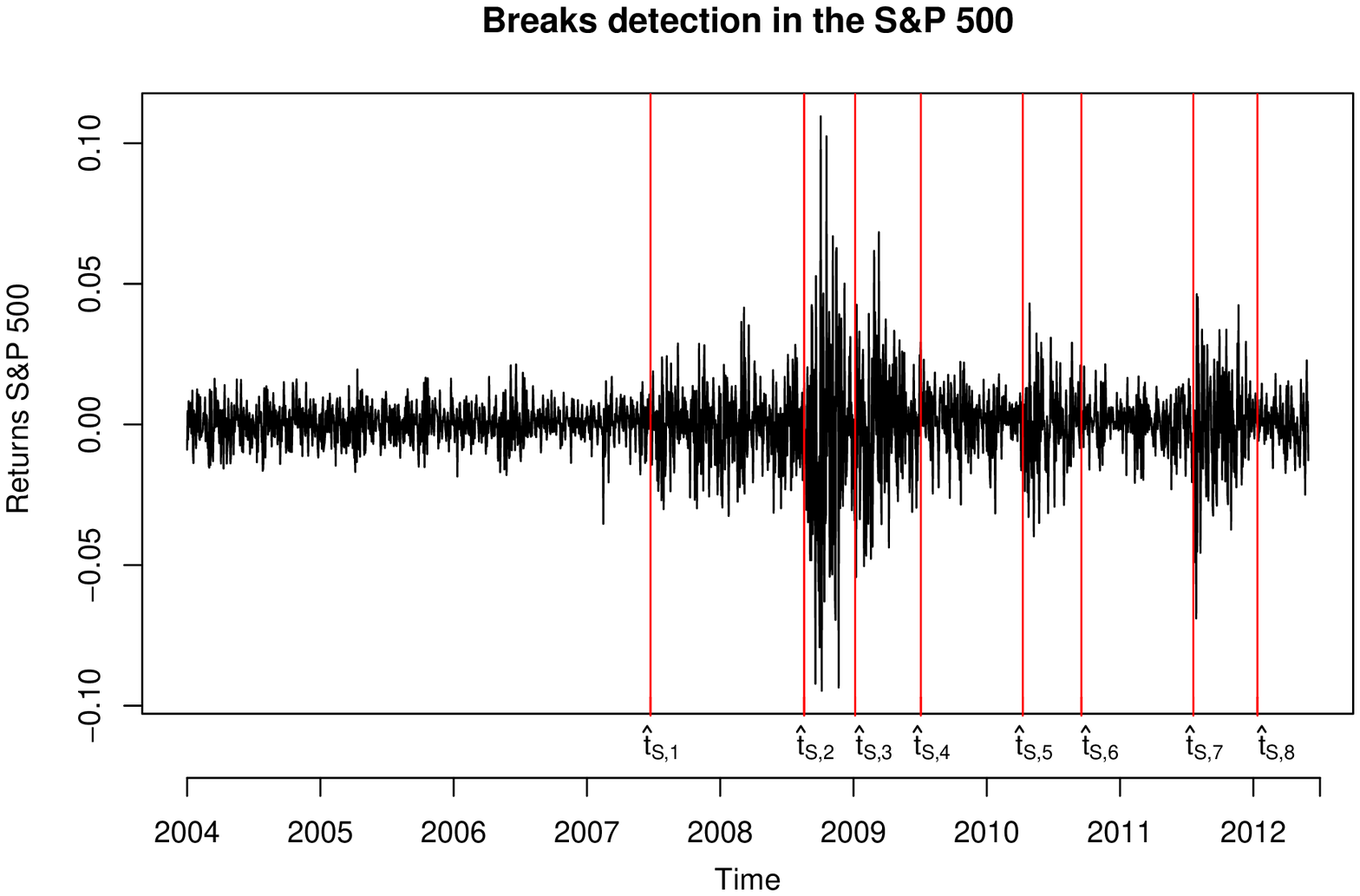 ,  width=15.1cm, height=9.5cm}\\
   \epsfig{file=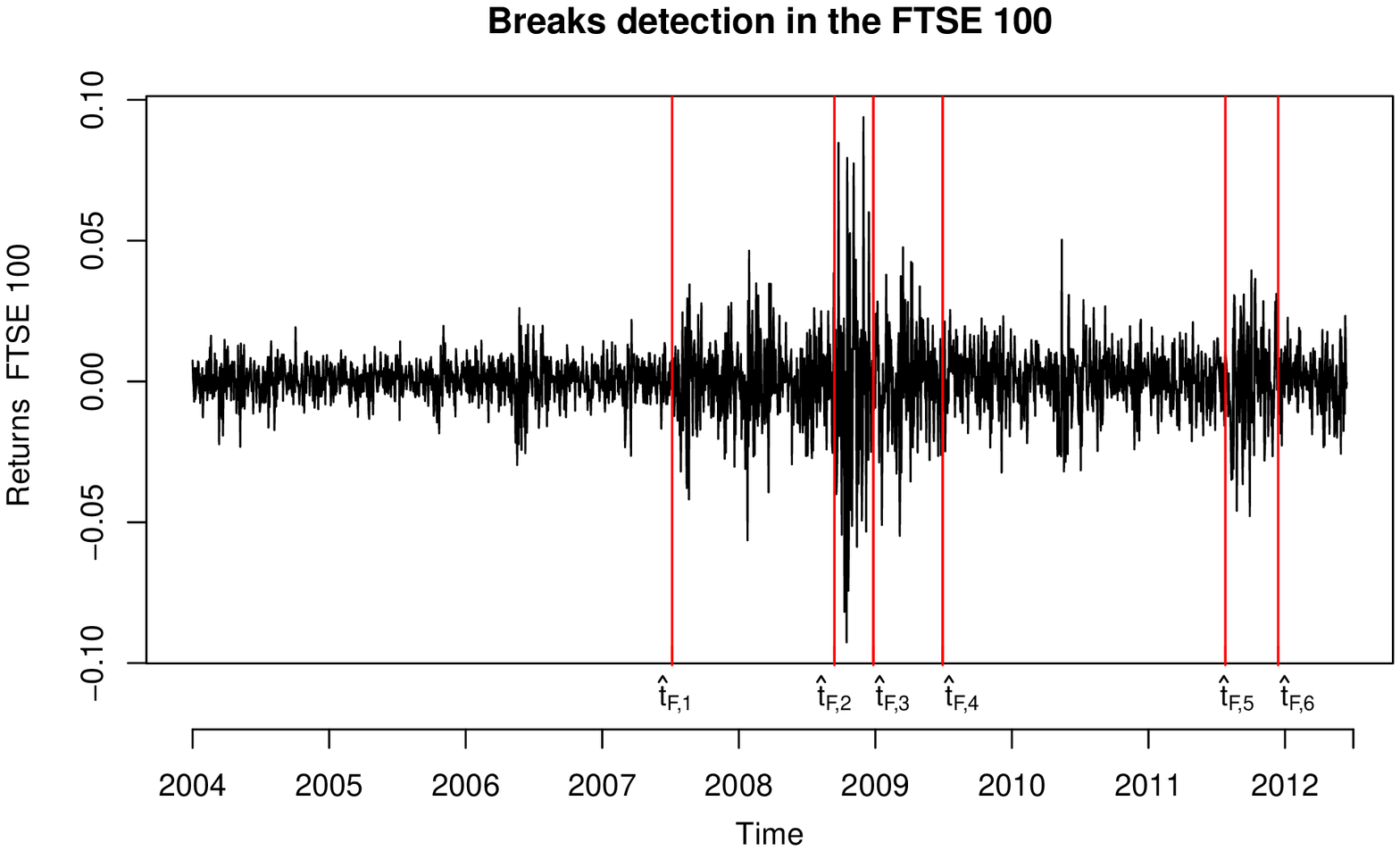, width=15.1cm, height=9.5cm}
   \end{center}
   \caption{ Break detection in the returns of S$\&$P 500 and FTSE 100 data using monitoring procedure based on $\widehat{C}_{k}$.
   The verticals  lines indicate the dates where breaks have been detected.}
   \label{fig6}
   \end{figure}

   We are going to see for S$\&$P 500 and FTSE 100 data how multiple changes can be monitored. For these series, monitoring starts at January 2, 2006.
   Figure \ref{fig3} represents the sequence $(\widehat{C}_k)$ for $k$ going from January 2, 2006 to December 31, 2008.
   According to the Figure \ref{fig3}, the monitoring stops at November 16, 2007 and September 4, 2007 for S$\&$P 500 and FTSE 100 data respectively.
   Retrospective procedure is applied to the series going from January 2, 2006 to November 16, 2007 (for S$\&$P 500 data) and the series
   going from January 2, 2006 to  September 4, 2007 (for FTSE 100 data). Breaks are detected at $\widehat{t}_{S,1} \simeq $ June 18, 2007 (for S$\&$P 500 data)
   and $\widehat{t}_{F,1} \simeq $ July 6, 2007  (for FTSE 100 data); see Figure \ref{fig6}. These breaks correspond to the beginning of the Subprime Crisis in US.\\
   After monitored the first change, we need to update the procedure. The new historical data are the series going from $\widehat{t}_{S,1}$
   to November 16, 2007 (for S$\&$P 500 data) and the series going from $\widehat{t}_{F,1}$  to September 4, 2007 (for FTSE 100 data).
   Therefore, monitoring continues at November 19, 2007 (for S$\&$P 500 data) and at September 5, 2007 (for FTSE 100 data).
   Figure \ref{fig4} shows the curve of the sequence $(\widehat{C}_{k})$.
   The monitoring stops at October 17, 2008 and November 10, 2008 for S$\&$P 500 and FTSE 100 data respectively.
   The retrospective test is applied and the break point estimation are $\widehat{t}_{S,2} \simeq $ August 14, 2008 and
   $\widehat{t}_{F,2} \simeq $ September 17, 2008 respectively for these two series;
   see Figure \ref{fig6}. These breaks correspond to the Lehman Brothers Bankruptcy which affects the worldwide financial system.\\
   \indent After that, the procedure is updated and monitoring continues at October 20, 2008 and November 11, 2008 for S$\&$P 500 and FTSE 100 data.
   Figure \ref{fig5} shows the sequence $(\widehat{C}_{k})$.
   The monitoring stopped at March 17, 2009 (S$\&$P 500) and March 9, 2009 (FTSE 100) and retrospective test
   detected change at $\widehat{t}_{S,3} \simeq $ January 5, 2009 (in S$\&$P 500) and $\widehat{t}_{F,3} \simeq $ December 29, 2008 (in FTSE 100 data); see Figure \ref{fig6}.
   These breaks correspond to the worldwide governments intervention to solve the financial crisis.\\
   The procedure continues until June 2012, other breaks are detected at $\widehat{t}_{S,4} \simeq $ 26 June 2009, $\widehat{t}_{S,5} \simeq $ 5 April 2010,
   $\widehat{t}_{S,6} \simeq $ 27 September 2010, $\widehat{t}_{S,7} \simeq $ 19 July 2011, $\widehat{t}_{S,8} \simeq $ 11 January 2012 (for S$\&$P 500 data) and
   $\widehat{t}_{F,4} \simeq $ 30 June 2009, $\widehat{t}_{F,5} \simeq $ 27 July 2011, $\widehat{t}_{F,6} \simeq $ 21 December 2011 (for FTSE 100 data).
   They are represented on Figure \ref{fig6}. These breaks correspond to the turmoils periods in the $2010-2012^+$ Greece and European debt crisis.

 ~ \\
   { \Large \textbf{Summary of the real-data applications} } \\
   Both monitoring procedure (based on detector $\widehat{C}_{k}$) and retrospective test have been applied to detect breaks
   in the Nikkei 225, S$\&$P 500 and FTSE 100  stock index. The following results are obtained :
   \begin{enumerate}
    \item For the Nikkei 225, from January 2, 1995 to October 19, 1998 ; break is detected at
    \begin{itemize}
        \item $\widehat{t}_N \simeq $ 17  September  1997 which correspond to the turmoil period  of the Asian financial crisis (1997-1998).
    \end{itemize}
    See also Figure \ref{fig2}.

    \item For the S$\&$P 500 and FTSE 100, from January 2, 2004 to June 11,
    2012 ;  break are detected at ($\widehat{t}_{S,i}$ and $\widehat{t}_{F,i}$ are referred to the breakpoint in the S$\&$P 500 and FTSE 100 respectively)
    \begin{itemize}
        \item $\widehat{t}_{S,1} \simeq $ 18 June 2007  and  $\widehat{t}_{F,1} \simeq$ 6 July 2007 which correspond to the beginning of the Subprime Crisis in US;
        \item $\widehat{t}_{S,2} \simeq $ 14 August 2008  and  $\widehat{t}_{F,2} \simeq$ 17 September 2008 which correspond to the Lehman Brothers Bankruptcy;
        \item  $\widehat{t}_{S,3} \simeq $ 5 January 2009    and  $\widehat{t}_{F,3} \simeq$  29 December 2008  which correspond
                 worldwide governments intervention to solve the financial crisis;
        \item $\widehat{t}_{S,4} \simeq $ 26 June 2009, $\widehat{t}_{S,5} \simeq $ 5 April 2010, $\widehat{t}_{S,6} \simeq $ 27 September 2010,
                $\widehat{t}_{S,7} \simeq $ 19 July 2011,  $\widehat{t}_{S,8} \simeq $ 11 January 2012 and
                $\widehat{t}_{F,4} \simeq $ 30 June 2009, $\widehat{t}_{F,5} \simeq $ 27 July 2011, $\widehat{t}_{F,6} \simeq $ 21 December 2011.
                These breaks indicates  the turmoils periods in the $2010-2012^+$ Greece and European debt crisis.
    \end{itemize}
   \end{enumerate}
  See also Figure \ref{fig6}.

 \newpage

\section{Proofs of the main results}
 Let us prove first some useful lemmas. In the sequel, for any $ x \in \R$, $[x]$ denotes the integer part of $x$.  Let $ (\psi_n)_n  $ and $ (r_n)_n  $ be sequences of random variables. Throughout this section, we use the notation
  $ \psi_n = o_P(r_n)  $ to mean :  for all $  \varepsilon  > 0, ~ P( |\psi_n| \geq \varepsilon|r_n| ) \rightarrow 0 $ as $n \rightarrow  \infty$.
  Write $ \psi_n = O_P(r_n)  $ to mean :  for all  $  \varepsilon > 0 $,  there exists $C>0$  such that  $ ~ P( |\psi_n| \geq C |r_n| )\leq \varepsilon  $
   for $n$ large enough.\\
 ~\\
  Recall that $(X_1,\cdots,X_n)$ is an observed trajectory of a process ${\cal M}_{\Z}(M_{\theta^*_0},f_{\theta^*_0})$. \\
 Let $k\geq n\geq 2$ and  $T_{1,n}=\{1,\cdots,n\}$,  $T_{\ell,k}=\{\ell,\ell+1,\cdots,k\}$ with $\ell \in \Pi_{n,k}=\{v_n,v_n+1,\cdots, k-v_n\}$, and
define
   $$ C_{k,\ell}:= \sqrt{n}\, \dfrac{k-\ell}{k}\big\|G^{-1/2}F\cdot\big(\widehat{\theta}(T_{\ell,k}) - \widehat{\theta}(T_{1,n})\big)\big\|,$$
   with $\widehat{\theta}$ defined in \eqref{theta}.
\begin{lem} \label{lem1}
 Under the assumptions of Theorem \ref{theo1},
 $$ \sup_{k>n}  ~\max_{ \ell \in \Pi_{n,k}} ~ \dfrac{1}{b((k-\ell)/n)}\, \big| \widehat{C}_{k,\ell} - C_{k,\ell}  \big|=o_P(1) ~ ~  \text{as} ~ n\rightarrow \infty.$$
 \end{lem}
 \begin{dem}
 For any $n\geq 1$, we have
  $$ \sup_{k>n}  ~\max_{ \ell \in \Pi_{n,k}}  ~ \dfrac{1}{b((k-\ell)/n)}\big| \widehat{C}_{k,\ell} - C_{k,\ell}  \big|
  =  \dfrac{1}{\inf_{s>0}b(s)} ~
  \sup_{k>n}  ~\max_{ \ell \in \Pi_{n,k}} ~ \big| \widehat{C}_{k,\ell} - C_{k,\ell}  \big| .$$
  Now, proceed similarly as in the proof of Lemma 3 of \cite{Kengne2011}.
 \end{dem}
\begin{lem} \label{lem2}
 Under the assumptions of Theorem \ref{theo1}
 \begin{multline*}
   \sup_{k>n}  ~\max_{ \ell \in \Pi_{n,k}}  ~ \dfrac{1}{b((k-\ell)/n)} \dfrac{\sqrt{n}}{k}
 \big\| (k-\ell)F\cdot\big( \widehat{\theta}(T_{\ell,k})- \widehat{\theta}(T_{1,n})\big)
  -2\big(  \dfrac{\partial}{\partial\theta}L(T_{\ell,k},\theta^*_0) - \dfrac{k-\ell}{n} \dfrac{\partial}{\partial\theta}L(T_{1,n},\theta^*_0) \big)  \big\| \\
  =o_P(1) ~ ~  \text{as} ~ n\rightarrow \infty.
    \end{multline*}
 \end{lem}
 \begin{dem}
  Let $k\geq n$  and $T\subset \{1,\cdots,k\}$.  By applying the Taylor expansion to the coordinates of $ \partial \widehat{L}(T,\cdot)/\partial \theta $,
 and using the fact that $ \partial \widehat{L}(T,\widehat{\theta}(T))/ \partial \theta = 0 $ we have
  $$ \dfrac{2}{\mbox{Card}(T)}\dfrac{\partial}{\partial\theta}\widehat{L}(T,\theta^*_0)=\widetilde{F}(T)\cdot(\widehat{\theta}(T) - \theta^*_0)
  ~ ~ \text{where}~ ~ \widetilde{F}(T)= -2\big( \dfrac{1}{\mbox{Card}(T)} \dfrac{\partial^2 \widehat{L}(T,  \widetilde{\theta}_{i}(T))}{ \partial \theta \partial
  \theta_i}\big)_{1\leq i \leq d} $$
  for some $\widetilde{\theta}_{i}(T)$ between  $\widehat{\theta}(T)$  and $\theta^*_0$.\\
  Hence for any $\ell \in \Pi_{n,k}$
  \begin{multline*}
  F \,(\widehat{\theta}(T_{\ell,k}) - \theta^*_0) =  \dfrac{2}{k-\ell}\dfrac{\partial}{\partial\theta}L(T_{\ell,k},\theta^*_0)
     + \big(F-  \widetilde{F}(T_{\ell,k}) \big)\big( \widehat{\theta}(T_{\ell,k}) - \theta^*_0  \big)\\
  + \dfrac{2}{k-\ell} \big(  \dfrac{\partial}{\partial\theta}\widehat{L}(T_{\ell,k},\theta^*_0) - \dfrac{\partial}{\partial\theta}L(T_{\ell,k},\theta^*_0)\big).
    \end{multline*}
    and
\begin{equation*}
  F \, (\widehat{\theta}(T_{1,n}) - \theta^*_0) =  \dfrac{2}{n}\dfrac{\partial}{\partial\theta}L(T_{1,n},\theta^*_0)
     + \big(F-  \widetilde{F}(T_{1,n})\big)\big( \widehat{\theta}(T_{1,n}) - \theta^*_0  \big)
  + \dfrac{2}{n} \big(  \dfrac{\partial}{\partial\theta}\widehat{L}(T_{1,n},\theta^*_0) - \dfrac{\partial}{\partial\theta}L(T_{1,n},\theta^*_0)\big).
    \end{equation*}
 Therefore, for any $\ell \in \Pi_{n,k}$
  \begin{multline}\label{eq1_lem2}
  \dfrac{\sqrt{n}}{k} \Big( (k-\ell)F\, \big( \widehat{\theta}(T_{\ell,k})- \widehat{\theta}(T_{1,n})\big)
  -2\big(  \dfrac{\partial}{\partial\theta}L(T_{\ell,k},\theta^*_0) - \dfrac{k-\ell}{n} \dfrac{\partial}{\partial\theta}L(T_{1,n},\theta^*_0) \big)  \Big) \\
  =  \sqrt{n} \dfrac{k-\ell}{k} \big(F-  \widetilde{F}(T_{\ell,k}) \big)\big( \widehat{\theta}(T_{\ell,k}) - \theta^*_0  \big)
  + 2\dfrac{\sqrt{n}}{k} \big(  \dfrac{\partial}{\partial\theta}\widehat{L}(T_{\ell,k},\theta^*_0) - \dfrac{\partial}{\partial\theta}L(T_{\ell,k},\theta^*_0)\big)\\
  - \sqrt{n} \dfrac{k-\ell}{k} \big(F-  \widetilde{F}(T_{1,n}) \big)\big( \widehat{\theta}(T_{1,n}) - \theta^*_0  \big)
  - 2\dfrac{k-\ell}{k} \dfrac{1}{\sqrt{n}} \big(  \dfrac{\partial}{\partial\theta}\widehat{L}(T_{1,n},\theta^*_0) - \dfrac{\partial}{\partial\theta}L(T_{1,n},\theta^*_0)\big).
    \end{multline}
 For $k>n$ and with some $\ell_k \in \Pi_{n,k}$, we have
 $$
 \underset{ \ell \in \Pi_{n,k}}{\mbox{max}} ~ \dfrac{1}{b((k-\ell)/n)} \sqrt{n} \dfrac{k-\ell}{k} \| F-  \widetilde{F}(T_{\ell,k})\|
\,  \| \widehat{\theta}(T_{\ell,k}) - \theta^*_0 \| \leq
\dfrac{1}{\underset{ s>0}{\mbox{Inf}}b(s) } \sqrt{k-\ell_k }\| F-
\widetilde{F}(T_{\ell_k,k}) \|
  \,  \|  \widehat{\theta}(T_{\ell_k,k}) - \theta^*_0 \|.
  $$
 According to \cite{Bardet2009} and \cite{Bardet2010}, $\| F-  \widetilde{F}(T_{\ell_k,k})\|=o_P(1)$ and
   $\|  \widehat{\theta}(T_{\ell_k,k}) - \theta^*_0 \| = O_P(1/\sqrt{k-\ell_k})$ as $k-\ell_k \rightarrow \infty$.
   Hence
    \begin{equation}\label{eq2_lem2}
   \underset{k>n} {\mbox{sup}} ~   \underset{ \ell \in \Pi_{n,k}}{\mbox{max}} ~ \dfrac{1}{b((k-\ell)/n)} \sqrt{n} \dfrac{k-\ell}{k}  \| F-  \widetilde{F}(T_{\ell,k}) \|
\,     \| \widehat{\theta}(T_{\ell,k}) - \theta^*_0 \| = o_P(1)  ~ ~
\text{as} ~ ~ n \rightarrow \infty .
   \end{equation}
  Similar arguments imply that
  \begin{equation}\label{eq3_lem2}
   \underset{k>n} {\mbox{sup}} ~ \underset{ \ell \in \Pi_{n,k}}{\mbox{max}} ~ \dfrac{1}{b((k-\ell)/n)} \sqrt{n} \dfrac{k-\ell}{k}  \| F-  \widetilde{F}(T_{1,n})\| \,
   \| \widehat{\theta}(T_{1,n}) - \theta^*_0 \| = o_P(1)  ~ ~ \text{as} ~ ~ n \rightarrow \infty .
    \end{equation}
 For $k>n$ and \text{for some} $\ell_k \in \Pi_{n,k}$,  we have
  \begin{multline*}
 \underset{ \ell \in \Pi_{n,k}}{\mbox{max}} ~ \dfrac{1}{b((k-\ell)/n)}  \dfrac{\sqrt{n}}{k}
  \|\dfrac{\partial}{\partial\theta}\widehat{L}(T_{\ell,k},\theta^*_0) - \dfrac{\partial}{\partial\theta}L(T_{\ell,k},\theta^*_0) \|
\\ \leq \dfrac{1}{\underset{ s>0}{\mbox{Inf}}b(s) } \dfrac{1}{\sqrt{k-\ell_k} } \|\dfrac{\partial}{\partial\theta}\widehat{L}(T_{\ell_k,k},\theta^*_0)
  - \dfrac{\partial}{\partial\theta}L(T_{\ell_k,k},\theta^*_0) \|.
    \end{multline*}
 According to \cite{Bardet2009},
  $\dfrac{1}{\sqrt{k-\ell_k} } \|\dfrac{\partial}{\partial\theta}\widehat{L}(T_{\ell_k,k},\cdot) - \dfrac{\partial}{\partial\theta}L(T_{\ell_k,k},\cdot) \|_{\Theta} = o_P(1)$
 as $k-\ell_k \rightarrow \infty$.
   Hence
   \begin{equation}\label{eq4_lem2}
   \underset{k>n} {\mbox{sup}} ~ \underset{ \ell \in \Pi_{n,k}}{\mbox{max}} ~ \dfrac{1}{b((k-\ell)/n)}  \dfrac{\sqrt{n}}{k}
  \|\dfrac{\partial}{\partial\theta}\widehat{L}(T_{\ell,k},\theta^*_0) - \dfrac{\partial}{\partial\theta}L(T_{\ell,k},\theta^*_0) \| = o_P(1)  ~ ~ \text{as} ~ ~ n \rightarrow \infty .
   \end{equation}
  Similar arguments show that
   \begin{equation}\label{eq5_lem2}
   \underset{k>n} {\mbox{sup}} ~ \underset{ \ell \in \Pi_{n,k}}{\mbox{max}} ~ \dfrac{1}{b((k-\ell)/n)}  \dfrac{k-\ell}{k} \dfrac{1}{\sqrt{n}}
  \|\dfrac{\partial}{\partial\theta}\widehat{L}(T_{1,n},\theta^*_0) - \dfrac{\partial}{\partial\theta}L(T_{1,n},\theta^*_0) \| = o_P(1)  ~ ~ \text{as} ~ ~ n \rightarrow \infty .
   \end{equation}
 Thus, Lemma \ref{lem2} follows from (\ref{eq1_lem2}), (\ref{eq2_lem2}), (\ref{eq3_lem2}), (\ref{eq4_lem2}) and (\ref{eq5_lem2}).
 \end{dem}
\begin{lem} \label{lem3}
 Under the assumptions of Theorem \ref{theo1}
 $$ \sup_{k>n}  ~\max_{ \ell \in \Pi_{n,k}}  ~ \dfrac{1}{ b((k-\ell)/n) }
 \sqrt{n} \dfrac{k-\ell}{k} \| F \cdot \big( \widehat{\theta}(T_{\ell,k})- \widehat{\theta}(T_{1,n})\big) \|
 \limiteloin \sup_{t>1}   \sup_{0<s<t}  ~ \dfrac{\| W_G(s)-sW_G(1) \|}{t~b(s)} $$
 where $W_G$ is a $d$-dimensional Gaussian centered process with covariance matrix $\E(W_G(s)W_G(\tau)')=\min(s,\tau)G$.
 \end{lem}
 \begin{dem}
  We are going to apply Lemma \ref{lem2} for specifying the asymptotic behaviour of $\widehat{\theta}(T_{\ell,k})- \widehat{\theta}(T_{1,n})$.

For $k>n$ and $\ell \in \Pi_{n,k} $, we have
  $$2 \dfrac{\sqrt{n}}{k} \big(  \dfrac{\partial}{\partial\theta}L(T_{\ell,k},\theta^*_0) - \dfrac{k-\ell}{n} \dfrac{\partial}{\partial\theta}L(T_{1,n},\theta^*_0) \big)
  = - \dfrac{n}{k} \dfrac{1}{\sqrt{n}} \big(  \sum_{i=\ell+1}^k \dfrac{\partial q_i(\theta^*_0)}{\partial\theta}
  -  \dfrac{k-\ell}{n} \sum_{i=1}^n \dfrac{\partial q_i(\theta^*_0)}{\partial\theta}  \big).$$
Now we are going to proceed in two steps.\\
 {\bf Step 1.}  Let $T>1$. We have
\begin{align*}
  &  \underset{n<k<nT} {\mbox{max}} ~\underset{ \ell \in \Pi_{n,k}}{\mbox{max}} ~ \dfrac{1}{ b((k-\ell)/n) }
     \dfrac{2\sqrt{n}}{k} \big \| \dfrac{\partial}{\partial\theta}L(T_{\ell,k},\theta^*_0) - \dfrac{k-\ell}{n} \dfrac{\partial}{\partial\theta}L(T_{1,n},\theta^*_0) \big \| \\
  &= \underset{n<k<nT} {\mbox{max}} ~\underset{ \ell \in \Pi_{n,k}}{\mbox{max}} ~ \dfrac{1}{ b((k-\ell)/n) }
  \dfrac{\sqrt{n}}{k} \big \|  \sum_{i=\ell}^{k} \dfrac{\partial  q_i(\theta^*_0)}{\partial\theta} - \dfrac{k-\ell}{n}\sum_{i=1}^{n} \dfrac{\partial  q_i(\theta^*_0)}{\partial\theta}   \big \| \\
  &= \max_{t\in \{1,1+\frac 1 n ,\cdots,T\}} \max _{s \in \{1-\frac {v_n} n,2-\frac {v_n} n,\cdots,t-\frac {v_n} n\}} \dfrac{1}{ b(([nt]-[ns])/n) } \, \frac n {[nt]}  \, \big \|     \dfrac{1}{\sqrt{n}}  \big( \sum_{i=[ns]+1}^{[nt]} \dfrac{\partial  q_i(\theta^*_0)}{\partial\theta}
  -  \dfrac{[nt]-[ns]}{n}\sum_{i=1}^{n} \dfrac{\partial  q_i(\theta^*_0)}{\partial\theta} \big) \big \|.
 \end{align*}
 Define the set $ S:=\{ (t,s)\in [1,T]\times [1,T]/ ~ s<t  \}  $. According to \cite{Bardet2009}, $\big( \dfrac{\partial q_i(\theta^*_0)}{\partial\theta} \big)_{t \in \Z}$  is a stationary ergodic martingale difference sequence with covariance matrix $G$. By Cramér-Wold device (see \cite{Billingsley1968} p. 206), it holds that
     $$  \dfrac{1}{\sqrt{n}}\sum_{i=[ns]+1}^{[nt]} \dfrac{\partial q_i(\theta^*_0)}{\partial\theta}
     ~ ~ \overset{ \mathcal{D}(S)}{\underset{n\rightarrow \infty} \longrightarrow} ~~  W_G(t-s).$$
 with $\overset{ \mathcal{D}(S)}{\underset{n\rightarrow \infty} \longrightarrow}$ means the weak convergence on the Skorohod space $\mathcal{D}(S)$.
 Hence
   $$
      \dfrac{1}{\sqrt{n}} \big( \sum_{i=[ns]+1}^{[nt]} \dfrac{\partial  q_i(\theta^*_0)}{\partial\theta}
  -  \dfrac{[nt]-[ns]}{n}\sum_{i=1}^{n} \dfrac{\partial  q_i(\theta^*_0)}{\partial\theta} \big) \\
   \overset{ \mathcal{D}(S)}{\underset{n\rightarrow \infty} \longrightarrow} ~~  W_G(t-s)-(t-s)W_G(1).
    $$
 Therefore
   \begin{eqnarray}
\nonumber && \underset{n<k<nT} {\mbox{max}} ~\underset{ \ell \in
\Pi_{n,k}}{\mbox{max}} ~ \dfrac{1}{ b((k-\ell)/n) }
  \dfrac{2\sqrt{n}}{k} \| \dfrac{\partial}{\partial\theta}L(T_{\ell,k},\theta^*_0) - \dfrac{k-\ell}{n} \dfrac{\partial}{\partial\theta}L(T_{1,n},\theta^*_0) \|    \\
 \nonumber && \hspace{5cm} \limiteloin \underset{1<t<T} {\mbox{sup}} ~  \underset{1<s<t}  {\mbox{sup}} ~ \dfrac{\| W_G(t-s)-(t-s)W_G(1) \|}{t~b(t-s)}\\
  && \hspace{5cm} \limiteloin \underset{1<t<T} {\mbox{sup}} ~  \underset{1<s<t}  {\mbox{sup}} ~ \dfrac{\| W_G(s)-s\, W_G(1) \|}{t~b(s)}.\label{eq1_lem3}
    \end{eqnarray}
  {\bf Step 2.} We will show that the limit distribution (as $n, T \rightarrow \infty$) of
  $$\underset{k>nT} {\mbox{sup}} ~\underset{ \ell \in \Pi_{n,k}}{\mbox{max}} ~ \dfrac{1}{ b((k-\ell)/n) }
  \dfrac{2\sqrt{n}}{k} \| \dfrac{\partial}{\partial\theta}L(T_{\ell,k},\theta^*_0) - \dfrac{k-\ell}{n} \dfrac{\partial}{\partial\theta}L(T_{1,n},\theta^*_0) \| $$
  exists and is equal to the limit distribution (as $ T \rightarrow \infty$) of
 $$\underset{t>T} {\mbox{sup}} ~  \underset{1<s<t}  {\mbox{sup}} ~ \dfrac{\| W_G(s)-s\, W_G(1) \|}{t~b(s)}.$$
  %
  Let $k>nT$. We have
  $$\underset{ \ell \in \Pi_{n,k}}{\mbox{max}} ~ \dfrac{1}{ b((k-\ell)/n) }  \dfrac{\sqrt{n}}{k} \| \dfrac{\partial}{\partial\theta}L(T_{\ell,k},\theta^*_0) \|
   \leq  \dfrac{1}{\underset{ s>0}{\mbox{Inf}}~b(s)} \dfrac{\sqrt{n}}{k} \|\sum_{i=\ell_k+1}^k \dfrac{\partial q_i(\theta^*_0)}{\partial\theta} \|
    ~ ~  \text{for some} ~ ~ \ell_k \in \Pi_{n,k} .$$
 It comes from the Hájek-Rényi-Chow inequality (see \cite{Chow1960}) that, for any $\varepsilon >0$
 $$  \lim_{T \rightarrow \infty} \lim_{n \rightarrow \infty }P\big(  \underset{k>nT} {\mbox{sup}}
 \dfrac{\sqrt{n}}{k} \|\sum_{i=\ell_k+1}^k \dfrac{\partial q_i(\theta^*_0)}{\partial\theta} \| > \varepsilon \big)=0 .$$
  Hence
  \begin{equation}\label{eq2_lem3}
  \underset{k>nT} {\mbox{sup}} ~\underset{ \ell \in \Pi_{n,k}}{\mbox{max}} ~ \dfrac{1}{ b((k-\ell)/n) }
  \dfrac{2\sqrt{n}}{k} \| \dfrac{\partial}{\partial\theta}L(T_{\ell,k},\theta^*_0) \| = o_P(1) ~ ~ \text{as} ~ ~ T,n \rightarrow \infty.
     \end{equation}
 Moreover, since the function $b(\cdot)$ is non-increasing,  for any $n,T>1$, we have:
   \begin{eqnarray}
  \underset{k>nT} {\mbox{sup}} ~\underset{ \ell \in \Pi_{n,k}}{\mbox{max}} ~ \dfrac{1}{ b((k-\ell)/n) } \dfrac{2\sqrt{n}}{k} \|  \dfrac{k-\ell}{n}\dfrac{\partial}{\partial\theta}L(T_{1,n},\theta^*_0) \|
 \nonumber   & = & \| \dfrac{1}{\sqrt{n}} \sum_{i=1}^n \dfrac{\partial q_i(\theta^*_0)}{\partial\theta}
                   \| \times \underset{k>nT} {\mbox{sup}} ~\underset{ \ell \in \Pi_{n,k}}{\mbox{max}} \dfrac{1}{ b((k-\ell)/n) }  \dfrac{k-\ell}{k}  \\
 \nonumber   & = & \| \dfrac{1}{\sqrt{n}} \sum_{i=1}^n \dfrac{\partial q_i(\theta^*_0)}{\partial\theta}
                   \| \times \underset{k>nT} {\mbox{sup}} ~ \dfrac{1}{ b((k-v_n)/n) }  \dfrac{k-v_n}{k}  \\
 \nonumber    &=&  \dfrac{1}{\underset{ s>0}{\mbox{Inf}}~b(s)} \| \dfrac{1}{\sqrt{n}} \sum_{i=1}^n \dfrac{\partial q_i(\theta^*_0)}{\partial\theta} \|    \\
  \label{eq3_lem3} &  & \hspace{-1cm} \limiteloin \dfrac{1}{\underset{ s>0}{\mbox{Inf}}~b(s)} \|W_G(1) \| ,
     \end{eqnarray}
     using again the Cramèr-Wold device.
 It comes from (\ref{eq2_lem3}) and (\ref{eq3_lem3}) that
   \begin{equation}\label{eq4_lem3}
  \underset{k>nT} {\mbox{sup}} ~\underset{ \ell \in \Pi_{n,k}}{\mbox{max}} ~ \dfrac{1}{ b((k-\ell)/n) }
  \dfrac{2\sqrt{n}}{k} \| \dfrac{\partial}{\partial\theta}L(T_{\ell,k},\theta^*_0) - \dfrac{k-\ell}{n} \dfrac{\partial}{\partial\theta}L(T_{1,n},\theta^*_0) \|
 ~ \overset{ \mathcal{D}}{\underset{T,n \rightarrow \infty} \longrightarrow} ~   \dfrac{1}{\underset{ s>0}{\mbox{Inf}}~b(s)} \|W_G(1) \|.
     \end{equation}
 Furthermore, since the coordinates of $W_G$ are Brownian motions, by the law of the iterated logarithm there exists $t_0>\exp(1)$ such as
 $$s>t_0 \Rightarrow \| W_G(s) \| \leq \sqrt{s}\log(s) ~ ~  \text{almost surely}.$$
 Thus,  for any $t>t_0$, we obtain almost surely
 $$\underset{1<s<t}{\mbox{sup}}\| W_G(s) \|  \leq \underset{1<s<t_0}{\mbox{sup}}\| W_G(s) \| + \sqrt{t}\log(t).$$
 Therefore, for  $T$ large enough, we have
 \begin{equation}\label{eq5_lem3}
  \underset{t>T}{\mbox{sup}} ~ \underset{1<s<t}{\mbox{sup}} \dfrac{\| W_G(s) \|}{t~b(s)} \leq \frac 1 {\inf _{s>0} b(s)} \, \big(  \dfrac{1}{T} ~ \underset{1<s<t_0}{\mbox{sup}}\| W_G(s) \|
   +  \underset{t>T}{\mbox{sup}} \dfrac{\log(t)}{\sqrt{t}} \big) ~ ~   \overset{ \text{a.s.}}{\underset{T \rightarrow \infty} \longrightarrow}  0 .
 \end{equation}
 Finally,  since $b(\cdot)$ is non-increasing,  for any $T>1$, we have
  \begin{equation}\label{eq6_lem3}
  \underset{t>T}{\mbox{sup}} ~ \underset{1<s<t}{\mbox{sup}} \dfrac{\| s W_G(1) \|}{t~b(s)} =\| W_G(1) \| ~\underset{t>T}{\mbox{sup}} \dfrac{1}{t} ~
  \underset{1<s<t}{\mbox{sup}}~ \dfrac{s}{b(t)} = \| W_G(1) \| ~\underset{t>T}{\mbox{sup}}~ \dfrac{1}{b(t)} =   \dfrac{1}{\underset{ s>0}{\mbox{Inf}}~b(s)} \|W_G(1) \| .
 \end{equation}
  It comes from (\ref{eq5_lem3}) and (\ref{eq6_lem3}) that the limit of \eqref{eq1_lem3} satisfies when $T\to \infty$,
   \begin{equation}\label{eq7_lem3}
  \underset{t>T}{\mbox{sup}} ~ \underset{1<s<t}{\mbox{sup}} \dfrac{\| W_G(s) - s W_G(1) \|}{t~b(s)}
 ~ \overset{ \mathcal{D}}{\underset{T \rightarrow \infty} \longrightarrow} ~   \dfrac{1}{\underset{ s>0}{\mbox{Inf}}b(s)} \|W_G(1) \|.
     \end{equation}
 From  {\bf Step 1} and {\bf Step 2} (the relations (\ref{eq1_lem3}), (\ref{eq4_lem3}) and  (\ref{eq7_lem3})), it comes that
 $$ \underset{k>nT} {\mbox{sup}} ~\underset{ \ell \in \Pi_{n,k}}{\mbox{max}} ~ \dfrac{1}{ b((k-\ell)/n) }
  \dfrac{2\sqrt{n}}{k} \| \dfrac{\partial}{\partial\theta}L(T_{\ell,k},\theta^*_0) - \dfrac{k-\ell}{n} \dfrac{\partial}{\partial\theta}L(T_{1,n},\theta^*_0) \|
 ~ \overset{ \mathcal{D}}{\underset{n \rightarrow \infty} \longrightarrow} ~   \underset{t>T}{\mbox{sup}} ~ \underset{1<s<t}{\mbox{sup}} \dfrac{\| W_G(s) - s W_G(1) \|}{t~b(s)}. $$
Hence, Lemma \ref{lem3} follows from  Lemma \ref{lem2}.
 \end{dem}

  \noindent {\bf Proof of Theorem \ref{theo1}} \\
   We know that
    \begin{align*}
 P\{\tau(n)<\infty\} & = P\Big\{ ~ \underset{k>n} {\mbox{sup}} ~ \underset{\ell \in \Pi_{n,k}} {\mbox{max}} \dfrac{\widehat{C}_{k,\ell}}{b((k-\ell)/n)}>1 \Big\} \\
   & =  P\Big\{~ \underset{k>n} {\mbox{sup}} ~\underset{ \ell \in \Pi_{n,k}}{\mbox{max}} ~ \dfrac{1}{ b((k-\ell)/n) }
             \sqrt{n} \dfrac{k-\ell}{k} \| \widehat{G}(T_{1,n})^{-1/2} \widehat{F}(T_{1,n}) \cdot \big( \widehat{\theta}(T_{\ell,k})- \widehat{\theta}(T_{1,n})\big) \| > 1 \Big\}.
   \end{align*}
  Since $ \widehat{G}(T_{1,n}) \limitepsn  G $ and $ \widehat{F}(T_{1,n}) \limitepsn  F $, it comes from Lemma \ref{lem1} and \ref{lem3} that
  \begin{multline*}\label{eq3_lem3}
 \underset{k>n} {\mbox{sup}} ~\underset{ \ell \in \Pi_{n,k}}{\mbox{max}} ~ \dfrac{1}{ b((k-\ell)/n) }
 \sqrt{n} \dfrac{k-\ell}{k} \| \widehat{G}(T_{1,n})^{-1/2} \widehat{F}(T_{1,n}) \cdot \big( \widehat{\theta}(T_{\ell,k})- \widehat{\theta}(T_{1,n})\big) \| \\
   \limiteloin   \underset{t>1} {\mbox{sup}} ~  \underset{1<s<t}  {\mbox{sup}} ~ \dfrac{\| G^{-1/2}( W_G(s)-sW_G(1)) \|}{t~b(s)} .
    \end{multline*}
 Since the covariance matrix of $ \{ W_G(s)~; ~ s \geq 0 \}$, is $\text{min}(s,\tau)G$, the covariance matrix of
 $ \{ G^{-1/2} W_G(s)~; ~ s \geq 0 \}$ is $\text{min}(s,\tau)I_d$ (where $I_d$ is the $d$-dimensional identity matrix).
 Hence Theorem \ref{theo1} follows. \Box
~\\

  \noindent {\bf Proof of Corollary \ref{cor1}} \\
  Since $b \equiv c$ a positive constant, it follows immediately from Theorem \ref{theo1} that
  $$\lim_{n \rightarrow \infty} P\{\tau(n)<\infty \}  = P\Big\{ \underset{t>1} {\mbox{sup}} ~ \underset{1<s<t} {\mbox{sup}} \frac{1}{t} \| W_d(s) - s W_d (1))\| >c \Big\}.$$
  Now, it suffices to show that
  $  \underset{t>1} {\mbox{sup}} ~ \underset{1<s<t} {\mbox{sup}} ~ \dfrac{1}{t} \| W_d(s) - s W_d (1))\| \overset{ \mathcal{D} } {=} U_d .$ \\
  For any $t>1$, we have
  $$ \underset{1<s<t} {\mbox{sup}} ~ \frac{1}{t}\, \|W_d(s) - s W_d (1))\| \overset{ \mathcal{D} } {=}
   \underset{1<s<t} {\mbox{sup}} ~ \frac{s}{t} \, \| W_d\big( \frac{s-1}{s} \big)\| = \underset{0<u< 1-1/t} {\mbox{sup}} ~ \frac{1}{t(1-u)} \, \| W_d(u)\| .
 $$
   Thus
    \begin{equation*}
   \underset{t>1} {\mbox{sup}} \underset{ ~ 1<s<t} {\mbox{sup}} ~ \frac{1}{t} \, \|W_d(s) - s W_d (1))\| \overset{ \mathcal{D} } {=}
     \underset{t>1} {\mbox{sup}} \underset{ ~ 0<u< 1-1/t} {\mbox{sup}} ~ \frac{1}{t(1-u)}\,  \| W_d(u)\|
     =  \underset{ 0< v<1} {\mbox{sup}} \underset{~ 0<u< v} {\mbox{sup}} ~ \frac{1-v}{1-u}\,  \| W_d(u)\|.
   \end{equation*}
 But,  $\| W_d(u)\|\overset{ \mathcal{D} } {=}  v^{1/2} \, \big \| W_d\big (\frac u v \big ) \big \|$. Therefore with $u=u'v$,
 \begin{equation*}
  \underset{ 0< v<1} {\mbox{sup}} \underset{~ 0<u< v} {\mbox{sup}} ~ \frac{1-v}{1-u}\,  \| W_d(u)\|= \underset{ 0< v<1} {\mbox{sup}} \underset{~ 0<u'<1 } {\mbox{sup}} ~ \frac{(1-v)v^{1/2}}{1-u'v}\,  \| W_d(u')\|.
   \end{equation*}
It remains to compute $ \displaystyle   \underset{ 0< v<1} {\mbox{sup}} \frac{(1-v)v^{1/2}}{1-u'v}$. Classical computations show that this supremum is obtained by $v=2 \, \big (3-u'+\sqrt{(9-u')(1-u')}\big )^{-1}$ and therefore
 \begin{multline}\label{eq2_cor1}
 \underset{t>1} {\mbox{sup}} \underset{ ~ 1<s<t} {\mbox{sup}} ~ \frac{1}{t} \, \|W_d(s) - s W_d (1))\| \overset{ \mathcal{D} } {=} \underset{~ 0<u'<1 } {\mbox{sup}} ~ f(u') \,  \| W_d(u')\| \\
 \quad \mbox{with}\quad f(u')=\frac {\sqrt{9-u'}+\sqrt{1-u'}}{\sqrt{9-u'}+3\sqrt{1-u'}}\, \Big ( \frac { 2}  {3-u'+\sqrt{(9-u')(1-u')}}\Big )^{1/2}.
   \end{multline}
     Hence,
   $$ \underset{t>1} {\mbox{sup}} ~ \underset{1<s<t} {\mbox{sup}} ~ \frac{1}{t} \| W_d(s) - s W_d (1))\|
    \overset{ \mathcal{D} } {=} U_d$$
     \Box

 \noindent {\bf Proof of Theorem \ref{theo2}} \\
   Denote $k_n= k^* + n^\delta$ for $\delta \in (1/2,1)$.
 For $n$ large enough, we have $v_n < n^\delta$ and thus $k_n - v_n = k^*+n^\delta -v_n \geq k^*$.
 Moreover, since $k^*>n$ then $k^* \in \Pi_{n,k} $ for $n$ large enough. \\
 In addition,  since $k^*=k^*(n)\geq n$ and  $\limsup_{n \rightarrow \infty} k^*(n)/n < \infty$,  there exists $c_0>1$ such that $k^*\leq c_0 n$ for $n$ large enough.
  Hence, according to assumption \textbf{B}, there exists a constant $c>0$ such that
   \begin{eqnarray}
   \nonumber \underset{ \ell \in \Pi_{n,k_n}}{\mbox{max}} \dfrac{\widehat{C}_{k_n,\ell}}{ b((k_n-\ell)/n)}
            & = & \underset{ \ell \in \Pi_{n,k_n}}{\mbox{max}} \dfrac{1}{ b((k_n-\ell)/n) }
            \sqrt{n} \dfrac{k_n-\ell}{k_n} \| \widehat{G}(T_{1,n})^{-1/2} \widehat{F}(T_{1,n}) \cdot \big( \widehat{\theta}_{k_n}(T_{\ell,k_n})- \widehat{\theta}(T_{1,n})\big) \\
 \nonumber & \geq  &  \dfrac{1}{ b((k_n-k^*)/n) }
            \sqrt{n} \dfrac{k_n-k^*}{k_n} \| \widehat{G}(T_{1,n})^{-1/2} \widehat{F}(T_{1,n}) \cdot \big( \widehat{\theta}_{k_n}(T_{k^*,k_n})- \widehat{\theta}(T_{1,n})\big) \\
   \nonumber   & \geq  &  c ~ \sqrt{n} ~ \dfrac{n^\delta}{k^* + n^\delta} \big\|\widehat{G}(T_{1,n})^{-1/2}\widehat{F}(T_{1,n})\big(\widehat{\theta}_{k_n}(T_{k^*,k_n}) - \widehat{\theta}(T_{1,n})\big)\big\|   \\
   \nonumber   & \geq  &     c  ~ \dfrac{n^{1/2+\delta}}{c_0 n + n^\delta} \big\|\widehat{G}(T_{1,n})^{-1/2}\widehat{F}(T_{1,n})\big(\widehat{\theta}_{k_n}(T_{k^*,k_n}) - \widehat{\theta}(T_{1,n})\big)\big\| \\
 \label{theo2_eq1}  & \geq  &   c  ~ \dfrac{n^{\delta - 1/2}}{(c_0 +1)} \big\|\widehat{G}(T_{1,n})^{-1/2}\widehat{F}(T_{1,n})\big(\widehat{\theta}_{k_n}(T_{k^*,k_n}) - \widehat{\theta}(T_{1,n})\big)\big\|.
   \end{eqnarray}

  According to \cite{Bardet2009} and \cite{Kengne2011}, $ \widehat{G}(T_{1,n}) \limitepsn G  $, ~   $ \widehat{F}(T_{1,n}) \limitepsn F $, ~
  $ \widehat{\theta}(T_{1,n}) \limitepsn \theta^*_0  $ ~ and ~ $ \widehat{\theta}_{k_n}(T_{k^*,k_n)}) \limitepsn \theta^*_1.$
  Since $G$ is symmetric positive definite, $F$ is invertible, $\theta^*_0 \neq \theta^*_1$ and $  \delta > 1/2 $ ,
    then (\ref{theo2_eq1}) implies that
    $$   \underset{\ell \in \Pi_{n,k_n}} {\mbox{max}} \dfrac{\widehat{C}_{k_n,\ell}}{b((k_n-\ell)/n)}  \limitepsn  \infty  .$$ \Box

\section*{Acknowledgements} The authors are grateful to the  anonymous referee for its report which helped to improve this paper.

\end{document}